\documentclass[12pt] {article}
\usepackage{amssymb}
\usepackage{amsmath,amsfonts}

\newtheorem {lemme} {{\bf Lemma}} [section]
\newtheorem {theoreme} {{\bf Theorem}} [section]
\newtheorem {proposition} {{\bf Proposition}} [section]

\newcommand{\tr}{\operatorname{tr}}
\newcommand{\Tr}{\operatorname{Tr}}
\newcommand{\N}{{\mathbb N}}
\newcommand{\R}{{\mathbb R}}
\newcommand{\C}{{\mathbb C}}
\newcommand{\vers}{\mathop{\longrightarrow}}
\numberwithin{equation}{section}
\title{ Strong asymptotic freeness for Wigner and Wishart  matrices}
\author{M. Capitaine\thanks
{Laboratoire de Statistique et Probabilit\'es, Universit\'e Paul Sabatier, 118 route de Narbonne, F-31062 Toulouse Cedex. E-mail: capitain@cict.fr }  \hspace{.1cm} and C. Donati-Martin\thanks
{Laboratoire de Probabilit\'es et Mod\`eles Al\'eatoires, Universit\'e Paris 6, Site Chevaleret, 13 rue Clisson, F-75013 Paris. E-mail: donati@ccr.jussieu.fr }}
\date{}
\begin{document}
\maketitle
\begin{abstract}
For each $n$ in $\N$, let $X_n = [(X_n)_{jk}]_{j,k=1}^n$ be  a
random  Hermitian matrix  such that  the $n^2$ random variables
$\sqrt{n}(X_n)_{ii}$, $\sqrt{2n} Re((X_n)_{ij})_{i<j}$, $\sqrt{2n}
Im((X_n)_{ij})_{i<j}$ are independent identically distributed with
common distribution $\mu  $ on $\R$. Let $X_n^{(1)}, \ldots,
X_n^{(r)}$ be $r$ independent copies of $X_n$ and $(x_1, \ldots
x_r)$ be a semicircular system in a $\cal{C}^{*}$-probability
space. Assuming that $\mu$ is symmetric and satisfies a Poincar\'e
inequality, we show that, almost everywhere, for any non commutative polynomial $p$
in $r$ variables,
\begin{equation} \label{VIfort}
\lim_{n \vers + \infty} || p(X_n^{(1)}, \ldots, X_n^{(r)}) || = ||
p(x_1, \ldots x_r)|| \ .
\end{equation}
We follow the method of \cite{HT} and \cite{S} which gave
(\ref{VIfort}) in the Gaussian (complex, real or symplectic) case.
We also get that (\ref{VIfort}) remains true when the $X_n^{(i)}$
are Wishart matrices while the $x_i$ are Marchenko-Pastur
distributed.
\end{abstract}

\noindent {\it Mathematics Subject Classification (2000)}: 15A52, 46L54, 60F99. \\
\noindent {\it Key words}:  Random matrices, free probability, asymptotic freeness.

\section{Introduction}
In the 90's, Voiculescu \cite{V}  introduced a random matrix model
for a free semi-circular system. He showed that if we take $r$
independent random matrices $(X_n^{(i)})_{i =1, \ldots r}$,
distributed as $GUE(n, \frac{1}{n})$, then, they are
asymptotically free, that is, for every non commutative polynomial
$p$ in $r$ variables,
\begin{equation} \label{V}
\mathbb{E}[\tr_n p(X_n^{(1)}, \ldots, X_n^{(r)})] \vers_{n \vers \infty} \tau(p(x_1, \ldots x_r))
\end{equation}
where $ tr_n$ stands  for the normalized trace on $M_n(\C)$ and $(x_1, \ldots x_r)$ is a free family of semicircular
variables in some non commutative probability space
$({\cal B}, \tau)$.  The result \eqref{V} holds true for a family of iid Wigner matrices and is proved by Dykema in \cite{D}.

\noindent
In a recent paper, Haagerup and Thorbj\o rnsen \cite{HT} proved a strong version of \eqref{V},
in the GUE case, namely a convergence for the operator norm:
\begin{equation} \label{Vfort}
\lim_{n \vers + \infty} || p(X_n^{(1)}, \ldots, X_n^{(r)}) || = || p(x_1, \ldots x_r)|| \ a.s.
\end{equation}
which led to the proof that $Ext(C^*_{red}(F_2))$ is not a group.

\noindent Schultz \cite{S} obtained the same result for Gaussian
random matrices in the real case (GOE) and in the simplectic case
(GSE). \noindent Our aim is to extend \eqref{Vfort} in the case of
an independent family
 of Wigner matrices on one hand  and in the case of Wishart matrices on the other hand .
Note that the special case $r=1$ gives the well known convergence
of the largest eigenvalue of $X_n^{(1)}$ to the right boundary of
the support of $x_1$ (see \cite{BY} for the Wigner case and
\cite{G} for the Wishart case; see also \cite{Ba} and the
references therein).

  Our approach is very similar to that of \cite{HT}
 and \cite{S}. Therefore, we will recall the main lines of their proofs. First, in proving \eqref{Vfort}, the minoration
$$ \liminf_{n \vers + \infty } || p(X_n^{(1)}, \ldots, X_n^{(r)}) || \geq  || p(x_1, \ldots x_r)|| \ a.s. $$
comes rather easily from an a.s. version of \eqref{V} (obtained in
\cite{Th} for the GUE case and proved in Section 6 of \cite{S} for
the GOE case) (see Lemma 7.2 in \cite{HT}). So, the main
difficulty is the proof of the reverse inequality:
\begin{equation}\label{limsup}
\limsup_{n \vers +\infty} || p(X_n^{(1)}, \ldots, X_n^{(r)}) ||
\leq  || p(x_1, \ldots x_r)|| \ a.s.
\end{equation}
In the following, we sketch the main steps in the proof of
(\ref{limsup}).

 \vspace{.3cm} \noindent
{\bf Step 1: A linearisation trick} (see \cite{HT}, Section 2 and Proposition 7.3) \\
In order to prove (\ref{limsup}), it is sufficient to prove:
\begin{lemme}   For all $m \in \N$, all self-adjoint matrices $a_0, \ldots, a_r$\footnote{By a density argument, we can also assume that the matrices
 $a_i$ are invertible.}
of size $m\times m$ and
all $\epsilon >0$,
\begin{equation} \label{spectre}
sp(a_0 \otimes 1_n + \sum_{i=1}^r a_i \otimes X_n^{(i)}(\omega)) \subset
sp(a_0 \otimes 1_{\cal B} + \sum_{i=1}^r a_i \otimes x_i) + ]-\epsilon, \epsilon[
\end{equation}
eventually, as $n \vers \infty$ a.e. in $\omega$. Here, $sp(T)$ denotes the spectrum of the operator $T$
and $1_n$ the identity matrix.
\end{lemme}
The analysis of the spectrum of $S_n := a_0 \otimes 1_n + \sum_{i=1}^r a_i \otimes X_n^{(i)}$ is done, using the Stieljes transform
\begin{equation} \label{Stieljes}
G_n( \lambda) = \mathbb{E}[(id_m \otimes \tr_n)[(\lambda \otimes
1_n - S_n)^{-1}]], \lambda \in M_m(\C), ~Im (\lambda)~
\mbox{positive definite}.
\end{equation}
The proof of \eqref{spectre} requires sharp estimates of the rate of convergence of $G_n(\lambda)$ to
$G( \lambda) := (id_m \otimes \tau)[(\lambda \otimes 1_{\cal B} - s)^{-1}]$ (of order $1/n^2$)
where $s = a_0 \otimes 1_{{\cal B}} + \sum_{p=1}^r a_p \otimes x_p$.

\vspace{.3cm} \noindent {\bf Step 2:} In the GUE case, Haagerup and
Thorbj\o rnsen \cite{HT} obtains the following estimate
\begin{equation} \label{mineq}
|| G_n( \lambda) - G(\lambda) || \leq \frac{C(\lambda)}{n^2}.
\end{equation}
 In the GOE, GSE cases, Schultz \cite{S} gets an extra term of
 order $1/n$, namely
\begin{equation} \label{GOEmi}
|| G_n( \lambda) - G(\lambda) - \frac{L(\lambda)}{n} || \leq
\frac{C(\lambda)}{n^2}
\end{equation}
for some functional $L$.

\vspace{.3cm} \noindent {\bf Step 3}  From the previous step, it
is shown in section 6 of \cite{HT} that
\begin{equation} \label{mve}
\mathbb{E}[(\tr_m \otimes \tr_n ) (\varphi(S_n))] = (\tr_m \otimes \tau)
(\varphi(s)) + O(\frac{1}{n^2})
\end{equation}
for $\varphi$ smooth with
compact support, and
\begin{equation} \label{5}
\mathbb{E}[(\tr_m \otimes \tr_n ) (\varphi(S_n))] =O(\frac{1}{n^2})
\end{equation}
for $\varphi$ smooth, constant outside a compact set and such that $supp(\varphi) \cap sp(s) = \emptyset$.
% with an estimation of $\frac{1}{n^4}$ for the variance. \\
In the GOE case (resp. GSE case), Schultz proved in section 5 of
\cite{S} that
\begin{equation} \label{mveL}
\mathbb{E}[(\tr_m \otimes \tr_n ) (\varphi(S_n))] = (\tr_m \otimes \tau)
(\varphi(s))+ \frac{1}{n}\Lambda(\varphi) + O(\frac{1}{n^2})
\end{equation}
where $\Lambda$ is a  distribution with compact support in $sp(s)$
with Stieljes transform
$$f(\lambda) = \tr_m(L(\lambda 1_m)), \lambda \in \C \setminus \R.$$
Therefore, \eqref{5} still holds for $\varphi$ with $supp(\varphi)
\cap sp(s) = \emptyset$.

\vspace{.3cm} \noindent {\bf Step 4}  (\ref{5}), combining with a
Gaussian variance estimate, yields  (by a standard application of
the Borel Cantelli lemma),
$$ (\tr_m \otimes \tr_n ) 1_F(S_n) = O(n^{-4/3})$$ for
$F = \{ t \in \R, d(t, sp(s)) \geq \epsilon \}$ which leads to \eqref{spectre}.

\vspace{.5cm} The main difficulties in the generalization of the
above to Wigner or Wishart matrices arise in step 2. Indeed, we don't have the gaussian integration by parts' formula anymore.
Our approach is inspired by the work of \cite{KKP} where they use a Taylor expansion (see Lemma \ref{IPP}) extending the gausian integration by parts' formula.
The remainder
of the proof can be completed essentially as in the GOE/GSE case.
Hence, in this paper, we shall focus on the obtention of such a
master inequality
$$ ||G_n(\lambda) - G(\lambda) - \frac{1}{n} L(\lambda) || = O(\frac{1}{n^2}) $$
in the case of a family of Hermitian matrices with symmetric iid
entries satisfying a Poincar\'e inequality, as well as in the case
of Wishart matrices;  we just give some hints when the
computations are similar to that of \cite{HT}, \cite{S}.

The paper is organized as follows. In section 2, we introduce
notations and preliminaries which will be of basic use later on.
In section 3, we describe the proof of (\ref{mve}) and (\ref{mveL}) proved
respectively in  \cite{HT} and \cite{S} in order to make clear the
validity of the method in our general framework we state in
section 4 (for the Wigner case) and section 5 (for the Wishart
case).

%%%%%%%%%%%%%%%%%%%%%%%%%
\section{Notations and preliminaries}
This section may contain some definitions already used in the introduction but we choose to
gather all the notations in this section for the reader's convenience.
To begin with, we introduce
some notations on the set of matrices.
\begin{itemize}
\item $M_p(\C)$ is the set of $p\times p$ matrices  with complex entries,  $M_p(\C)_{sa}$ the subset of self-adjoint elements of $M_p(\C)$ and $1_p$ the identity matrix. In the following, we shall consider two sets of matrices with $p=m$ ($m$ fixed) and $p=n$ with $n\vers \infty$.
\item $\Tr_p$ denotes the trace and $\tr_p = \frac{1}{p} \Tr_p$ the normalized trace on $M_p(\C)$.
\item $|| . ||$ denotes the operator norm on $M_p(\C)$ and $||M||_2 = (\Tr_p (M^*M))^{1/2}$ the Hilbert-Schmidt norm.
 \item  Let $(E_{ij})_{i,j =1}^n$ be the canonical basis of $M_n(\C)$ and define a basis of the real vector space of the self-adjoint matrices $M_n(\C)_{sa}$ by:
 \begin{eqnarray*}
 e_{jj} &= &E_{jj}, 1\leq j \leq n \\
 e_{jk}& = &\frac{1}{\sqrt{2}}(E_{jk} + E_{kj}), 1 \leq j < k \leq n \\
 f_{jk} &=& \frac{i}{\sqrt{2}}(E_{jk} - E_{kj}) , 1 \leq j < k \leq n
 \end{eqnarray*}
 \item[-] For a matrix $M$ in $M_m(\C) \otimes M_n(\C)$, we denote by
 $$M_{ij} : = (id_m \otimes \Tr_n)(M (1_m \otimes E_{ji})) \in M_m(\C), \ 1\leq i,j \leq n$$
 and
 $$ _{\alpha, \beta}\!M := (\Tr_m \otimes id_n)(M (\hat{E}_{\beta, \alpha}\otimes 1_n)) \in M_n(\C), \  1 \leq \alpha, \beta \leq m $$
  where $(\hat{E}_{\alpha, \beta})$ is the canonical basis of $M_m(\C)$.
  \end{itemize}
We now define  our matrix model and the random variables of interest.
\begin{itemize}
\item[-] $(X^{(1)}_n, \ldots, X^{(r)}_n)_{i = 1, \ldots r}$ is a
set of iid random matrices in $M_n(\C)_{sa}$, whose distribution
will be specified later (matrices in GUE or GOE in section 3, Wigner matrices in Section 4, Wishart
matrices in section 5).
 \item[-] For a given family $a_0, \ldots a_r$ in $M_m(\C)_{sa}$, we define the random variable $S_n$ with values in $M_m(\C) \otimes M_n(\C)$ by:
 \begin{equation} \label{defSn}
 S_n = a_0 \otimes 1_n + \sum_{p=1}^r a_p \otimes X^{(p)}_n
 \end{equation}
 and $s \in M_m(\C) \otimes {\cal B}$ by
  \begin{equation} \label{defS}
 s = a_0 \otimes 1_{{\cal B}} + \sum_{p=1}^r a_p \otimes x_p
 \end{equation}
 where the $(x_i)_{i = 1, \ldots r}$ is a free family of self-adjoint operators in a $C^*$ probability space
 $({\cal B},  \tau)$ with a faithful state $\tau$,  whose distribution will be specified in the different cases
(semi-circular in sections 3 and 4 or distributed as the Marchenko-Pastur distribution in section 5).
 \item[-] For any matrix $\lambda$ in ${\cal O}$ where
 $${\cal O} : = \{\lambda \in M_m(\C) | Im(\lambda) \mbox{ is positive definite} \}
 %=  \{\lambda \in M_m(\C) | \lambda_{min}(Im\lambda)>0 \}
 , $$
  we define the $M_m(\C)$ valued rv:
 \begin{equation} \label{defSt}
 H_n(\lambda) = (id_m \otimes \tr_n) [ (\lambda \otimes 1_n - S_n)^{-1}],
 \end{equation}
 \begin{equation} \label{espSt}
G_n(\lambda) = \mathbb{E}[ H_n(\lambda) ]
\end{equation}
and
\begin{equation} \label{defStlimite}
 G(\lambda) = (id_m \otimes \tau) [ (\lambda \otimes 1_{{\cal B}} - s)^{-1}].
 \end{equation}
 For $\lambda \in \C \setminus \R$, we also define
 $$g_n(\lambda) = \tr_m(G_n(\lambda 1_m))$$
 and $$g(\lambda) = \tr_m(G(\lambda 1_m)).$$
%\item[-] It is important to know the dependence in $\lambda$ in the estimations \eqref{mineq} or (\ref{GOEmi}) in order
% to pass from Step 2 to Step 3.

 \end{itemize}
 We end this preliminary by recalling some properties of $G(\lambda)$ and of the resolvent $(\lambda \otimes 1_n - S_n)^{-1} $ of the  matrix $S_n$.
 First, one can easily see that for any $\lambda$ and $\lambda^{'}$ in $M_m(\C)$ such that $Im(\lambda)$ and $Im(\lambda^{'})$ are positive
definite,
 \begin{equation}\label{difStieljes}(\lambda \otimes 1_{{\cal B}} - s)^{-1}-(\lambda^{'} \otimes 1_{{\cal B}} -
 s)^{-1}=(\lambda \otimes 1_{{\cal B}} - s)^{-1}(\lambda^{'} -
 \lambda)(\lambda^{'} \otimes 1_{{\cal B}} - s)^{-1}.
 \end{equation}
\begin{lemme} \label{G}
Let $\lambda$ in $M_m(\C)$ such that $Im(\lambda)$ is positive
definite. Then
 \begin{equation}\label{normeG}
 \Vert (\lambda \otimes 1_{{\cal B}} - s)^{-1}\Vert \leq ||
 Im(\lambda)^{-1}\Vert \mbox{~~and~~~}
  || G(\lambda ) || \leq || Im(\lambda)^{-1}
  ||.\end{equation}
Moreover, $G(\lambda )$ is invertible and
\begin{equation}\label{Gmoins}
  || G(\lambda )^{-1} || \leq( ||\lambda ||+ ||s|| )^2 || Im(\lambda)^{-1}
  ||.\end{equation}
\end{lemme}
\noindent We refer the reader to section 5 of \cite{HT} for a
proof of (\ref{Gmoins}).
 \begin{lemme} \label{lem2}
 Let $\lambda$ in $M_m(\C)$ such that $Im(\lambda)$ is positive definite, then
 \begin{equation}\label{norme}
  || (\lambda \otimes 1_n - S_n)^{-1} || \leq || Im(\lambda)^{-1}
  ||,\end{equation}
\begin{equation}  \label{2.9jk}
\forall 1 \leq k,l \leq n, || (\lambda \otimes 1_n - S_n)^{-1}_{kl} || \leq || Im(\lambda)^{-1} ||,
\end{equation}
  and for $p \geq 2$,
 \begin{equation} \label{2.10}
 \frac{1}{n} \sum_{k,l =1}^n ||(\lambda \otimes 1_n - S_n)^{-1}_{kl} ||^p \leq C_m ||Im(\lambda)^{-1} ||^p
 \end{equation}
 where, in the first inequality,  $||. ||$ denotes the operator norm in $M_m(\C)\otimes M_n(\C)$ ( in $M_m(\C)$ in the others)  and $C_m$ a constant depending only on $m$.

 \vspace{.2cm}
 \noindent
 For a Hermitian matrix $M$, the derivative w.r.t $M$ of the resolvent $R(z) = (z-M)^{-1}$ satisfies:
 \begin{equation} \label{resolvente}
 R'_M(z) . A = R(z) A R(z)   \mbox{ for all  Hermitian matrix $A$}. \end{equation}
 \end{lemme}
 {\bf Sketch of Proof:} We just mention the proof of \eqref{2.10}. From \eqref{2.9jk}, it's enough to consider the case $p=2$. \\
 Let us denote $G^{(n)} =  (\lambda \otimes 1_n - S_n)^{-1} \in M_m(\C)\otimes M_n(\C)$.
 Since the operator norm is smaller than the Hilbert-Schmidt norm,
 \begin{eqnarray*}
  \frac{1}{n} \sum_{k,l =1}^n ||G^{(n)}_{kl} ||^2 & \leq &  \frac{1}{n} \sum_{k,l =1}^n \sum_{\alpha, \beta = 1}^m |  _{\alpha, \beta}G^{(n)}_{kl} |^2 \\
  &=& \frac{1}{n} \Tr_{nm} (G^{(n)} (G^{(n)})^*) \\
  & \leq& m || G^{(n)} (G^{(n)})^* || \leq m  || Im(\lambda)^{-1} ||^2.   \end{eqnarray*}
  where the last inequality follows from \eqref{norme}. $\Box$
 %%%%%%%%%%%%%%%%%

 \vspace{.3cm}
 \noindent
  In the sequel, we shall denote by $P_k$ any polynomial of degree $k$
 whose coefficients are positive and by $C$ or $K$ any constant; $P_k$, $C$ or $K$   can depend on the
 $a_l$, $l=1, \ldots,r$,
 and may vary from line to line.
%%%%%%%%%%%%%%%%%%%%%%%%%%
\section{Main ideas in the proofs of (\ref{mve}) and (\ref{mveL}) from \cite{HT}  and \cite{S} }
\subsection{Estimate of $\Vert G_n-G\Vert$ in \cite{HT}}
Let us recall the main ideas of \cite{HT} in the estimation of
$\Vert G_n(\lambda) -G(\lambda)\Vert$. In lemma 5.4 of \cite{HT},
Haagerup and Thorbj\o rnsen observe in one hand that the
matrix-valued Stieljes transform of $s$ satisfies, for  any
$\lambda$ in ${\cal O}$,
\begin{equation}\label{eqS}
 \sum_{i=1}^r a_i G( \lambda)a_i  \ + \ (a_0 - \lambda)
+G (\lambda)^{-1}=0. \end{equation} In the other hand, using the
Gaussian integration by parts formula, they establish the analogue
of (\ref{eqS}) satisfied by $H_n(\lambda)$ (``Master equation",
Lemma 2 \cite{HT}):
\begin{equation}\label{MEG}
\mathbb{E}\left[ \sum_{p=1}^r   a_p H_n(\lambda) a_p  H_n(\lambda)
+ (a_0 - \lambda ) H_n(\lambda) + 1_m \right] = 0.\end{equation}
Then, using the Gaussian Poincar\'e inequality to get an estimate
of the variance of $H_n(\lambda)$, they deduce from (\ref{MEG})
the ``Master inequality" (Lemma 3 in \cite{HT}):
\begin{equation}\label{MI}
 \Vert \sum_{i=1}^r a_i G_n( \lambda)a_i G_n( \lambda)  \ + \ (a_0 - \lambda)
G (\lambda) + 1_m\Vert \leq \frac{C}{n^2}\Vert Im(\lambda)^{-1}
\Vert^4. \end{equation} Moreover, the authors prove that
$G_n(\lambda)$ is invertible for any $\lambda$ in ${\cal O}$ and
they give an upper bound of the norm of its inverse (Proposition
5.2 \cite{HT}) $$\Vert G_n(\lambda)^{-1}\Vert\leq (\Vert \lambda
\Vert +K)^2 \Vert Im(\lambda)^{-1} \Vert. $$ Hence, they deduce
from (\ref{MI}) that, for any $\lambda$ in ${\cal O}$,
\begin{equation}\label{MIdif}
 \Vert a_0+ \sum_{i=1}^r a_i G_n( \lambda)a_i   \ + \  G_n(\lambda)^{-1} - \lambda
\Vert \leq f_n(\Vert Im(\lambda)^{-1}\Vert,\Vert \lambda \Vert),
\end{equation}
\noindent where
$$    f_n(\Vert Im(\lambda)^{-1}\Vert,\Vert \lambda \Vert)    =\frac{C}{n^2} (\Vert \lambda
\Vert +K)^2  \Vert Im(\lambda)^{-1} \Vert^5.$$ Further, they set
$$\Lambda_n(\lambda) =  a_0+ \sum_{i=1}^r a_i G_n( \lambda)a_i   \ + \
G_n(\lambda)^{-1}$$ for any $\lambda$ in ${\cal O}$. (\ref{MIdif})
can be rewritten
\begin{equation}\label{dif}
 \Vert \Lambda_n(\lambda) - \lambda
\Vert \leq f_n(\Vert Im(\lambda)^{-1}\Vert,\Vert \lambda \Vert).
\end{equation}
The authors define
$${\cal O}_n^{'}=\{\lambda \in {\cal O}, f_n(\Vert Im(\lambda)^{-1}\Vert, \Vert \lambda
\Vert)< \frac{\epsilon(\lambda)}{2}\}$$ where
$$\epsilon(\lambda):=\frac{1}{\Vert Im(\lambda)^{-1}\Vert}.$$
(\ref{dif}) implies that, for any $\lambda$ in ${\cal O}_n^{'}$,
\begin{equation}\label{binf}
Im \Lambda_n(\lambda) \geq \frac{1}{2\Vert
Im(\lambda)^{-1}\Vert}1_m
\end{equation}
\noindent and that in particular $\Lambda_n(\lambda)$ belongs to
${\cal O}$ (see Lemma 5.5 \cite{HT}). Consequently, applying
(\ref{eqS}), they get that, for any $\lambda$ in ${\cal O}_n^{'}$,
\begin{equation}\label{idn}
a_0+ \sum_{i=1}^r a_i G( \Lambda_n(\lambda))a_i   \ + \
G(\Lambda_n(\lambda))^{-1}=a_0+ \sum_{i=1}^r a_i G_n( \lambda)a_i
\ + \ G_n(\lambda)^{-1}.
\end{equation}
In proof of (b) Proposition 5.6, Haagerup and Thorbj\o rnsen show
that (\ref{idn}) implies that
\begin{equation}\label{eg}G_n(\lambda)=G(\Lambda_n(\lambda))\end{equation}
\noindent for any $\lambda$ in ${\cal O}_n^{''}:=\{\lambda \in
{\cal O}_n^{'}, \epsilon(\lambda)> \sqrt{2\sum_{i=1}^r \Vert a_i
\Vert^2} \}.$ Using that $t \mapsto f_n(t^{-1},t) t^{-1}$ is a
continuous strictly decreasing function from $]0; +\infty[$ onto
$]0; +\infty[$, they show in proof of (a) Proposition 5.6
\cite{HT} that $ {\cal O}_n^{'}$ is an open connected subset of
$M_m(\C)$. Thus, by the principle of uniqueness of analytic
continuation, (\ref{eg}) still holds for any $\lambda$ in ${\cal
O}_n^{'}$. Thus, for any $\lambda$ in ${\cal O}_n^{'}$, they get
that \begin{eqnarray*} \Vert G_n(\lambda) -G(\lambda)\Vert& \leq &
\Vert G(\Lambda_n(\lambda)) -G(\lambda)\Vert\\
&\leq & \Vert Im(\Lambda_n(\lambda))^{-1}\Vert \Vert \lambda -
\Lambda_n(\lambda) \Vert \Vert Im(\lambda)^{-1}\Vert\\
& \leq & 2 f_n(\Vert Im(\lambda)^{-1}\Vert,\Vert \lambda \Vert)
\Vert Im(\lambda)^{-1}\Vert^2
\end{eqnarray*}
\noindent where the last inequality comes from (\ref{dif}),
(\ref{binf}). Now, if $\lambda$ belongs to ${\cal O} \setminus
{\cal O}_n^{'}$, they note that
\begin{eqnarray*} \Vert G_n(\lambda) -G(\lambda)\Vert& \leq &
  2\Vert Im(\lambda)^{-1}\Vert\\
& \leq & 4 f_n(\Vert Im(\lambda)^{-1}\Vert,\Vert \lambda \Vert)
\Vert Im(\lambda)^{-1}\Vert^2
\end{eqnarray*}
\noindent since
$$\frac{1}{2} \leq f_n(\Vert Im(\lambda)^{-1}\Vert,\Vert \lambda \Vert)
\Vert Im(\lambda)^{-1}\Vert.$$ Finally, for any $\lambda $ in
$\cal{O}$,
\begin{eqnarray} \Vert G_n(\lambda) -G(\lambda)\Vert
& \leq & 4 f_n(\Vert Im(\lambda)^{-1}\Vert,\Vert \lambda \Vert)
\Vert Im(\lambda)^{-1}\Vert^2 \nonumber\\
&=&\frac{C}{n^2}(\Vert \lambda \Vert +K)^2 \Vert Im(\lambda)^{-1}
\Vert^7.\label{estimHT}
\end{eqnarray}
%%%%%%%%%%%%%%%%%%%%%
\subsection{Estimate of $\Vert G_n-G-\frac{1}{n}L\Vert$ in \cite{S}}
In the GOE case, a term of order $1/n$ appears in the Master
equation so that the estimate of $\Vert G_n(\lambda)
-G(\lambda)\Vert$ Schultz makes by sticking to the previous proof
of \cite{HT} is of order $1/n$. Nevertheless, a further study (we
will describe in our general framework in section 4) gives her the
sharper estimate
\begin{equation}\label{estimLS}\Vert G_n(\lambda)-G(\lambda)-\frac{1}{n}L(\lambda)\Vert
\leq \frac{1}{n^2}(\Vert \lambda \Vert +K)^8 P_{13}(\Vert
Im(\lambda)^{-1}\Vert ) \end{equation} for any $\lambda$ such that $Im
\lambda$ positive definite or negative definite.
\subsection{From Step 2 to Step 3}
>From the previous estimates (\ref{estimHT}) and (\ref{estimLS}),
Haagerup, Thorbj\o rnsen and Schultz immediately get that, for any
$\lambda$ in $\C \setminus \R$,
\begin{equation}\label{estimgdif}
\vert r_n(\lambda)\vert \leq \frac{1}{n^2}(\Vert \lambda \Vert
+K)^{\alpha} P_{k}(\vert Im(\lambda)^{-1} \vert) \end{equation}
\noindent where

-in the GUE case \cite{HT}
$$r_n(\lambda) = g_n(\lambda)-g(\lambda) \ ,\ \alpha =2\,\ k=7.$$
-in the GOE case \cite{S}
$$r_n(\lambda) = g_n(\lambda)-g(\lambda)-\frac{1}{n}tr_m(L(\lambda 1_m)) \ ,\ \alpha =8\,\ k=13.$$
Since $S_n$ and $s$ are selfadjoint, by the spectral
theory, there exist unique probability measures $\mu_n$ and $\mu$
on $\R$ such that
$$\int \varphi d\mu_n = \mathbb{E}[(\tr_m \otimes \tr_n ) (\varphi(S_n))] $$
$$\int \varphi d\mu = (\tr_m \otimes \tau) (\varphi(s)).$$
\noindent $g_n$ and $g$ are the Stieljes transforms of $\mu_n$ and
$\mu$. Moreover, in Lemma 5.5 in \cite{S}, Schultz proves by using
a characterisation theorem of Tillmann that $l(\lambda) :=
tr_m(L(\lambda1_m))$ is the Stieljes transform of a distribution
$\Lambda$ with compact support in $sp(s)$. Hence, using the
inverse Stieljes tranform, Haagerup, Thorbjornsen and Schultz get
respectively that, for any $\varphi $ in ${\cal C}_c^{\infty}(\R,
\R)$,

\noindent - in \cite{HT} \begin{equation}\label{StHT} \int \varphi
d\mu_n-\int \varphi d\mu= -\frac{1}{\pi} \lim_{y \rightarrow 0^+}
Im \int_{\R} \varphi(x) r_n(x+iy)dx.\end{equation} - in \cite{S}
\begin{equation}\label{StS} \int \varphi d\mu_n-\int \varphi
d\mu -\frac{\Lambda(\varphi)}{n}= -\frac{1}{\pi} \lim_{y
\rightarrow 0^+} Im \int_{\R} \varphi(x)
r_n(x+iy)dx.\end{equation} Hence, the remainder of the two proofs
(in \cite{HT} and \cite{S}) deals with the estimation of $$
\limsup_{y \rightarrow 0^+} \vert \int_{\R} \varphi(x)
h(x+iy)dx\vert$$ \noindent  where $h$ is an analytic function on
$\C \setminus \R$ which satisfies
\begin{equation}\label{nestimgdif}
\vert h(\lambda)\vert \leq (\Vert \lambda \Vert +K)^{\alpha}
P_{k}(\vert Im(\lambda)^{-1} \vert).\end{equation} In \cite{HT}
section 6, Haagerup and Thorbj\o rnsen introduce a very clever
family of functions $\{I_p(\lambda), p \geq 1\}$ defined by
$$I_p(\lambda) =\frac{1}{(p-1)!}\int_0^{+\infty} h(\lambda+t)
t^{p-1} \exp(-t) dt.$$ They note that
$$I_1(\lambda) - I_1^{'}(\lambda)=h(\lambda)$$
$$I_p(\lambda) - I_p^{'}(\lambda)=I_{p-1}(\lambda),\ p \geq 2,
$$
so that for any $\varphi $ in ${\cal C}_c^{\infty}(\R, \R)$ and
$y>0$,
$$\int_{\R} \varphi(x)
h(x+iy)dx=\int_{\R}(1+D)^p \varphi(x) I_p(x+iy)dx.$$ Now, they
choose $p=k+1$ where $k$ is the degree of the polynomial in the
right hand side of (\ref{nestimgdif}) (that is $p=8$ in \cite{HT}
and $p=14$ in \cite{S}) and estimate $I_{k+1}(\lambda)$ for $Im
\lambda >0$. Using (\ref{nestimgdif}), it is not difficult to see
that
$$\lim_{r \rightarrow +\infty } \int_{[r,r+ir]}\frac{1}{k!} h(\lambda+z)
z^{k} \exp(-z) dz=0.$$ Thus, by Cauchy's integral theorem, the
authors get
\begin{eqnarray*}
I_{k+1}(\lambda)&=&\lim_{r \rightarrow +\infty }
\int_{[0,r+ir]}\frac{1}{k!} h(\lambda+z) z^{k} \exp(-z) dz\\
&=&  \int_{0}^{+\infty}\frac{1}{k!} h(\lambda+(1+i)t)
(1+i)^{k+1}t^k \exp(-(1+i)t) dt.
\end{eqnarray*}
Plugging in (\ref{nestimgdif}), one gets for any $\lambda$ such
that $Im \lambda >0$,
\begin{eqnarray*}
\vert I_{k+1}(\lambda)\vert &\leq&
\frac{2^{\frac{k+1}{2}}}{k!}\int_{0}^{+\infty}(\vert \lambda\vert
+ \sqrt{2} t +K)^{\alpha} P_k(\vert Im \lambda +t\vert^{-1})t^k
\exp(-t)dt\\
&\leq& \frac{2^{\frac{k+1}{2}}}{k!}\int_{0}^{+\infty}(\vert
\lambda\vert + \sqrt{2} t +K)^{\alpha} P_k( t^{-1})t^k
\exp(-t)dt\\
&\leq& \int_{0}^{+\infty}(\vert \lambda\vert + \sqrt{2} t
+K)^{\alpha} Q(t) \exp(-t)dt
\end{eqnarray*}
where $Q(t) = \frac{2^{\frac{k+1}{2}}}{k!}P_k(t^{-1})t^k$ is a
polynomial.

\noindent It follows by dominated convergence
$$\limsup_{y \rightarrow 0^+} \vert \int_{\R} \varphi(x)
h(x+iy)dx\vert \leq \int_{\R}\int_{0}^{+\infty}\vert (1+D)^p
\varphi(x)\vert (\vert x\vert + \sqrt{2} t
+K)^{\alpha}Q(t)\exp(-t)dt dx < +\infty.$$ \noindent Dealing with
$h(\lambda) =n^2 r_n(\lambda)$ one gets

\begin{equation} \label{majlimsup} \limsup_{y \rightarrow 0^+} \vert \int_{\R} \varphi(x)
r_n(x+iy)dx\vert \leq \frac{C}{n^2}. \end{equation}

\noindent Combining (\ref{majlimsup}) with respectively
(\ref{StHT}) and (\ref{StS}) , one gets respectively (\ref{mve}) and (\ref{mveL}).

 \section{The iid case}
 We consider a Hermitian matrix $X_n =  [(X_n)_{jk}]_{j,k=1}^n$  of size n for which the $n^2$ rv $((X_n)_{ii})$,
$(\sqrt{2} Re((X_n)_{ij})_{i<j}$, $(\sqrt{2} Im((X_n)_{ij})_{i<j}$
are independent identically distributed with common distribution
$\mu / \sqrt{n}$ where $\mu$ is a symmetric distribution with
variance 1 on $\R$ which satisfies a Poincar\'e inequality (see
section 4.2). We call $X_n$ a Wigner matrix with distribution
$\mu$. Let $X_n^{(1)}, \ldots, X_n^{(r)}$ be $r$ independent
copies of $X_n$. \noindent We present our main technical tool (see
\cite{KKP}):
 \begin{lemme} \label{lem1}
Let $\xi$ be a real-valued rv such that  $\mathbb{E}(\vert \xi
\vert^{p+2})<\infty$. Let  $\phi$ be a function from  $\R$ to $\C$
such that the first $p+1$ derivatives are  continuous and bounded.
Then,
\begin{equation}\label{IPP}\mathbb{E}(\xi \phi(\xi)) = \sum_{a=0}^p
\frac{\kappa_{a+1}}{a!}\mathbb{E}(\phi^{(a)}(\xi)) +
\epsilon\end{equation} where  $\kappa_{a}$ are the cumulants of
$\xi$, $\vert \epsilon \vert \leq C \sup_t \vert
\phi^{(p+1)}(t)\vert \mathbb{E}(\vert \xi \vert^{p+2})$, $C$
depends on $p$ only.
\end{lemme}
In the following, we shall apply this identity with a function
$\phi(\xi)$ given by the Stieljes transform of a random matrix. It
follows from the Lemma \ref{lem2} and (\ref{resolvente}) above
that the conditions of Lemma \ref{lem1} (bounded derivatives) are
fulfilled.

 \subsection{The master equation}
 Note that since $\mu$ satisfies a Poincar\'e inequality, we have
 $\int \vert x \vert^q d \mu(x) < + \infty $ for any $q$ in $\N$ (see
 Corollary 3.2 and Proposition 1.10  in \cite{L}). Note also that,
since  $\mu$ is symmetric, any odd cumulant of $\mu$ vanishes.
 \begin{theoreme} \label{ME} With the previous notations,
 \begin{equation} \label{mastereq}
 \mathbb{E} \left[ \sum_{i=1}^r a_i H_n( \lambda)a_i H_n(\lambda) \ + \ (a_0 - \lambda) H_n (\lambda) +1_m \right]
 + \frac{1}{n} R_n(\lambda) +\epsilon_n= 0
 \end{equation}
 where $\Vert \epsilon_n \Vert \leq \frac{P_6 (\Vert Im(\lambda)^{-1} \Vert)}{n^{2}}$ and
 $R_n(\lambda)$ denotes the quantity
 $$ \frac{\kappa_4}{2}  \mathbb{E} \left[\sum_{p=1}^r
 \frac{1}{n^2} \sum_{k,l =1}^n a_p (\lambda \otimes 1_n - S_n)^{-1}_{kk} a_p (\lambda \otimes 1_n - S_n)^{-1}_{ll}
 a_p (\lambda \otimes 1_n - S_n)^{-1}_{kk}
   a_p  (\lambda \otimes 1_n - S_n)^{-1}_{ll} \right]
 $$
 where $\kappa_4$ is the fourth cumulant of the distribution $\mu$.
 Note that \begin{equation}\label{majR} \Vert R_n(\lambda) \Vert \leq P_4 (\Vert Im(\lambda)^{-1}
 \Vert).\end{equation}
 \end{theoreme}
 {\bf Proof:}
 We shall apply formula \eqref{IPP} to the ${\cal M}_m(\C)$-valued function $\phi(\xi) = (\lambda \otimes 1_n - S_n)^{-1}_{ij}$ for $1 \leq i,j \leq n$ and $\xi$ is one
 of the variable $(X_n^{(p)})_{kk}$,
$\sqrt{2} Re((X_n^{(p)})_{kl})$, $\sqrt{2} Im((X_n^{(p)})_{kl})$ for $1 \leq k<l \leq n$ and $p \leq r$. \\
We notice that
\begin{eqnarray*}
\frac{\partial \phi}{\partial Re((X_n^{(p)})_{kk})}  &=& \phi'_{X_n^{(p)}} . e_{kk}, \ 1 \leq k \leq n \\
\frac{\partial \phi}{\partial \sqrt{2} Re((X_n^{(p)})_{kl})}  &=& \phi'_{X_n^{(p)}} . e_{kl}, \ 1 \leq k< l \leq n\\
\frac{\partial \phi}{\partial \sqrt{2} Im((X_n^{(p)})_{kl})}  &=& \phi'_{X_n^{(p)}} . f_{kl} , \ 1 \leq k < l \leq n\\
\sqrt{2} Re((X_n^{(p)})_{kl}) &=& \Tr_n(X_n^{(p)} e_{kl}), \ 1 \leq k <l \leq n \\
\sqrt{2} Im((X_n^{(p)})_{kl}) &=& \Tr_n(X_n^{(p)} f_{kl}), \ 1 \leq k <l \leq n \\
(X_n^{(p)})_{kk} &=& \Tr_n(X_n^{(p)} e_{kk}), \ 1 \leq k \leq n.
\end{eqnarray*}

\noindent Let $1 \leq p\leq r$, $1\leq k\leq l \leq n$  be fixed.
For simplicity, we write $\phi'$, $\phi''$, $\phi'''$ for the
first derivatives of $\phi$ with respect to $\sqrt{2}
Re((X_n^{(p)})_{kl})$. Then, according to (\ref{resolvente}),
\begin{eqnarray*}
\phi' &=& \left[ (\lambda \otimes 1_n - S_n)^{-1} a_p \otimes e_{kl} (\lambda \otimes 1_n - S_n)^{-1} \right]_{ij} \\
\phi''&=&2 \left[ (\lambda \otimes 1_n - S_n)^{-1} a_p \otimes e_{kl} (\lambda \otimes 1_n - S_n)^{-1}  a_p \otimes e_{kl} (\lambda \otimes 1_n - S_n)^{-1}\right]_{ij} \\
\phi'''&=&6 \left[ (\lambda \otimes 1_n - S_n)^{-1} a_p \otimes e_{kl} (\lambda \otimes 1_n - S_n)^{-1}  a_p \otimes e_{kl} (\lambda \otimes 1_n - S_n)^{-1} \right.\\
&&\left.  \quad a_p \otimes e_{kl} (\lambda \otimes 1_n - S_n)^{-1}\right]_{ij}
\end{eqnarray*}
Writing \eqref{IPP} in this setting gives
\begin{equation} \label{1}
\mathbb{E}[ \Tr_n(X_n^{(p)} e_{kl}) (\lambda \otimes 1_n -
S_n)^{-1}_{ij} ] = \frac{1}{n} \mathbb{E}[ \phi'] +
\frac{\kappa_4}{6n^2} \mathbb{E}[\phi''']  + O(n^{-3})
\end{equation}
where the $O(n^{-3})$ means the norm of this term is smaller than
 $\frac{C\Vert a_p\Vert^5 \Vert Im(\lambda)^{-1} \Vert^6}{n^{3}}$.
Multiplying by $n$ gives the equation, denoted by $A_{ij}^{kl}(p)$:
\begin{equation} \label{2}
n \mathbb{E}[ \Tr_n(X_n^{(p)} e_{kl}) (\lambda \otimes 1_n -
S_n)^{-1}_{ij} ] = \mathbb{E}[ \phi'] + \frac{\kappa_4}{6n}
\mathbb{E}[\phi'''] + O(n^{-2})
\end{equation}
with the analogous equations with $f_{pq}$ (denoted by $B_{ij}^{kl}(p)$) and $e_{pp}$. \\

\noindent
Recall how we can obtain the master equation in the gaussian case (GUE case) from  \eqref{2} which reads in this case:
\begin{equation} \label{3}
n \mathbb{E}[ \Tr_n(X_n^{(p)} e_{kl}) (\lambda \otimes 1_n -
S_n)^{-1}_{ij} ] = \mathbb{E}\left[ (\lambda \otimes 1_n -
S_n)^{-1} a_p \otimes e_{kl} (\lambda \otimes 1_n - S_n)^{-1}
\right]_{ij}.
\end{equation}
By a linear combination with the analogous equation with $f_{kl}$, we have:
$$n \mathbb{E}[ \Tr_n(X_n^{(p)} E_{kl}) (\lambda \otimes 1_n - S_n)^{-1}_{ij} ] = \mathbb{E}\left[ (\lambda \otimes 1_n - S_n)^{-1} a_p \otimes E_{kl} (\lambda \otimes 1_n - S_n)^{-1} \right]_{ij}$$
for all $1\leq k,l \leq n$. \\
Now, take in the above formula $i =k$, $j=l$ and consider $\frac{1}{n^2} \sum_{k,l}$, we then obtain:
\begin{equation} \label{GUE}
\frac{1}{n} \sum_{k,l} \mathbb{E}[ (X_n^{(p)})_{lk} (\lambda
\otimes 1_n - S_n)^{-1}_{kl} ] =  \frac{1}{n^2} \sum_{k,l}
\mathbb{E}\left[ (\lambda \otimes 1_n - S_n)^{-1}_{kk} a_p
(\lambda \otimes 1_n - S_n)^{-1}_{ll} \right]
\end{equation}
that is
\begin{equation} \label{4}
 \mathbb{E}[ id_m \otimes \tr_n ((1_m \otimes X_n^{(p)}) (\lambda \otimes 1_n - S_n)^{-1}) ] = \mathbb{E}\left[ H_n(\lambda) a_p  H_n(\lambda)  \right].
\end{equation}
Now, from the above equation,
\begin{eqnarray*}
 \mathbb{E}\left[ \sum_{p=1}^r  a_p H_n(\lambda) a_p  H_n(\lambda)  \right] &=&
 \sum_{p=1}^r \mathbb{E}[ id_m \otimes \tr_n ((a_p \otimes 1_n)(1_m \otimes X_n^{(p)}) (\lambda \otimes 1_n - S_n)^{-1}) ] \\
 &=& \sum_{p=1}^r \mathbb{E}[ id_m \otimes \tr_n ((a_p  \otimes X_n^{(p)}) (\lambda \otimes 1_n - S_n)^{-1}) ] \\
 &=& \mathbb{E}[ id_m \otimes \tr_n (S_n - a_0 \otimes 1_n) (\lambda \otimes 1_n - S_n)^{-1}) ] \\
 &=& -1_m  + (\lambda - a_0 ) \mathbb{E}[H_n(\lambda)]
 \end{eqnarray*}
 implying the master formula in the GUE case:
 $$ \mathbb{E}\left[ \sum_{p=1}^r   a_p H_n(\lambda) a_p  H_n(\lambda) + (a_0 - \lambda ) H_n(\lambda) + 1_m \right] = 0.$$
 Keeping in mind these computations, we now study the terms coming from
  third derivatives. \\
We thus consider $A(p) = \frac{1}{n^2} \sum_{k,l} A_{kl}^{kl}(p)$ (resp. $B(p)$) and study all the contributions of the different terms.

\vspace{.3cm}
\noindent

\noindent
{\bf Study of the third derivative} \\
Writing as before the terms appearing in $\phi'''$, we can see that all the terms except one contains at least two $G_{kl}$ and then, according to  Lemma \ref{lem2}, these terms will give a contribution in $O(n^{-2})$ in $A(p)$. The only term to be considered  is:
$$\frac{1}{n^2} \sum_{k,l}  \mathbb{E}((\lambda \otimes 1_n - S_n)^{-1}_{kk} a_p  (\lambda \otimes 1_n - S_n)^{-1}_{ll} a_p (\lambda \otimes 1_n - S_n)^{-1}_{kk} a_p (\lambda \otimes 1_n - S_n)^{-1}_{ll}).$$

\vspace{.3cm}
\noindent
Now, using the same linear combination giving \eqref{GUE} in the GUE case, we obtain that the corrective term
 of order $1/n$ appearing in the iid case is:
\begin{eqnarray*}
&& \frac{1}{n}\{\frac{\kappa_4}{2} \mathbb{E}\left[ \frac{1}{n^2}
\sum_{k,l =1}^n (\lambda \otimes 1_n - S_n)^{-1}_{kk} a_p (\lambda
\otimes 1_n - S_n)^{-1}_{ll} a_p (\lambda \otimes 1_n -
S_n)^{-1}_{kk}
 \right. \\
&& \qquad  \qquad   \qquad  \qquad \left.   a_p  (\lambda \otimes
1_n - S_n)^{-1}_{ll} \right]\}.
 \end{eqnarray*}
 The proof of the Theorem is complete. $\Box$
 %%%%%%%
 %%%%%%%%%%
 \subsection{Variance estimate}
 We assume that $\mu$ satisfies a  Poincar\'e
inequality: there exists a positive constant $C$ such that for any
${\cal C}^1$ function $f: \R \rightarrow \C$  such that $f$ and
$f'$ are in $L^2(\mu)$,
$$\mathbf{V}(f)\leq C\int \vert f' \vert^2 d \mu ,$$
\noindent with $\mathbf{V}(f) = \mathbb{E}(\vert
f-\mathbb{E}(f)\vert^2)$.  We refer the reader to \cite{B}  for a
characterization of the measures on $\R$ which satisfy a
Poincar\'e inequality (see also \cite{Tou}). For example, $\mu (dx) = \exp(-|x|^\alpha) dx$ with $\alpha \geq 1$ satisfies the Poincar\'e inequality.

\noindent For any matrices $A_1, \ldots, A_r$, define $\Vert (A_1,
\ldots, A_r) \Vert_e^2 := \sum_{j=1}^r \Vert A_j\Vert_2^2$. Let
$\Psi: (M_n(\C)_{sa})^r \rightarrow \mathbb{R}^{rn^2}$ be the
canonical isomorphism introduced in
 Remark 3.4 in \cite{HT}.

 \begin{lemme}\label{variance}

For any  function $f: \R^{rn^2} \rightarrow \C$ be a ${\cal
C}^1$-function  such that $f$ and the gradient $\nabla(f )$ are both polynomially
bounded,
\begin{equation}\label{Poincare}\mathbf{V}{\left[ f\circ \Psi (X_n^{(1)}, \ldots,
X_n^{(r)})\right]} \leq \frac{C}{n}\mathbb{E}\{\Vert \nabla\left[f
\circ \Psi(X^{(1)}, \ldots, X^{(r)})\right] \Vert_e^2
\}.\end{equation}
\end{lemme}
{\bf Proof}
 $\mu^{(n)}:=\mu/\sqrt{n}$ satisfies the Poincar\'e inequality
$$\int \vert g-\int g d \mu^{(n)}\vert^2 d\mu^{(n)}\leq \frac{ C}{n}\int \vert g'\vert ^2 d \mu^{(n)}.$$
(\ref{Poincare}) readily follows by the tensorisation property of
the Poincar\'e inequality.
%%%%%%%
 %%%%%%%%%%
\subsection{Master inequality}

We follow  the lines of the proof of Theorem 4.5 in \cite{HT}.
Using
the master equality (\ref{mastereq}), we easily get\\

$ \Vert
\sum_{i=1}^r a_i G_n( \lambda)a_i G_n(\lambda) \ + \ (a_0 -
\lambda) G_n (\lambda) +1_m  +
 \frac{1}{n} R_n(\lambda) + \epsilon_n \Vert $
 $$\leq
\Vert \sum_{i=1}^r a_i^2\Vert \mathbb{E}\{\Vert
H_n(\lambda)-\mathbb{E}( H_n(\lambda))\Vert^2\}.$$

 \noindent Thanks to (\ref{Poincare}), the following
of the proof of Theorem 4.5 in \cite{HT} still holds and we
similarly get
$$ \mathbb{E}\{\Vert H_n(\lambda)-E(H_n(\lambda))\Vert^2\}\leq\frac{C m^3}{n^2}
\Vert \sum_{i=1}^r a_i^2\Vert \Vert (Im(\lambda))^{-1}\Vert^4$$
and therefore
\begin{equation}\label{MastIn}
\Vert  \sum_{i=1}^r a_i G_n( \lambda)a_i G_n(\lambda) \ + \ (a_0 -
\lambda) G_n (\lambda) +1_m  + \frac{1}{n} R_n(\lambda) \Vert
\leq\frac{P_6(\Vert (Im(\lambda))^{-1}\Vert)}{n^2} .
\end{equation}

\subsection{Estimation of $\Vert G_n-G\Vert$}
In the Gaussian case, Haagerup and Thorbj\o rnsen in \cite{HT} and
Schultz in \cite{S} prove that $G_n(\lambda)$ is invertible for
any $\lambda$ such that $Im \lambda$ is positive definite. In our
more general case, we are going to use the master inequality
(\ref{MastIn}) in order to prove that, for any $\lambda$ in some
subset of ${\cal O}$, $G_n(\lambda)$ is invertible and to get an
upper bound of $\Vert G_n(\lambda)^{-1}\Vert$.

\noindent Set $$B_n(\lambda)=  \sum_{i=1}^r a_i G_n( \lambda)a_i
 \ + \ (a_0 - \lambda).$$
% We have $$Im B_n(\lambda)=  \sum_{i=1}^r a_i Im G_n( \lambda)a_i
 %\ + \ - Im \lambda.$$
 %Since $Im G_n( \lambda)$ is negative (see...) and $- Im \lambda$
 %is definite negative, we deduce that
 %$Im B_n(\lambda)$ is definite negative and thus $ B_n(\lambda)$
%is invertible by lemma 3.1 in \cite{HT}.
Now, from the master
inequality (\ref{MastIn}) and (\ref{majR}), we get
\begin{equation}\label{MastIn2}
\Vert  \sum_{i=1}^r a_i G_n( \lambda)a_i G_n(\lambda) \ + \ (a_0 -
\lambda) G_n (\lambda) +1_m   \Vert \leq \frac{P_6(\Vert
(Im(\lambda))^{-1}\Vert)}{n}.
\end{equation}
\noindent that is $$\Vert B_n(\lambda)G_n (\lambda)+1_m   \Vert
\leq \frac{P_6(\Vert (Im(\lambda))^{-1}\Vert)}{n}.$$ \noindent
Hence, for any $\lambda$ such that $\frac{P_6(\Vert
(Im(\lambda))^{-1}\Vert)}{n}<\frac{1}{2}$, $B_n(\lambda)G_n
(\lambda)$ is invertible with
$$\Vert [B_n(\lambda)G_n(\lambda)]^{-1}\Vert \leq 2.$$
Thus, for such a $\lambda$, $G_n(\lambda)$ is also obviously
invertible with
\begin{equation}\label{majGinvmat}
\Vert G_n(\lambda)^{-1}\Vert \leq 2 \Vert B_n(\lambda)\Vert  \leq
2(\Vert a_0\Vert + \Vert \lambda \Vert +  \sum_{i=1}^r \Vert a_i
\Vert^2 \Vert (Im(\lambda))^{-1}\Vert)
\end{equation}
%\noindent where we used that
%$$\Vert  \sum_{i=1}^r a_i G_n( \lambda)a_i \Vert \leq \Vert \sum_{i=1}^r a_i^2
%\Vert \Vert G_n(\lambda)\Vert \leq \Vert \sum_{i=1}^r a_i^2 \Vert
%\Vert (Im(\lambda))^{-1}\Vert$$ which follows from the fact that
%$v \mapsto  \sum_{i=1}^r a_i va_i$ is completely positive and then
%atteins it norm at the unit of $M_m(C)$.
Now, from the inequality (\ref{MastIn2}) and using
(\ref{majGinvmat}), we get readily that for any $\lambda$ in
${\cal O}$ such that $\frac{P_6(\Vert
(Im(\lambda))^{-1}\Vert)}{n}<\frac{1}{2}$,
\begin{equation}\label{MastIn3}
\Vert  \sum_{i=1}^r a_i G_n( \lambda)a_i  \ + \ (a_0 - \lambda)
+G_n (\lambda)^{-1}    \Vert \leq \frac{P_6(\Vert
(Im(\lambda))^{-1}\Vert)}{n}2(\Vert a_0\Vert + \Vert \lambda \Vert
+  \sum_{i=1}^r \Vert a_i \Vert^2 \Vert (Im(\lambda))^{-1}\Vert).
\end{equation}
\noindent Define $${\cal O}^{'}_n=\left\{\lambda \in {\cal O},
\frac{P_6(\Vert (Im(\lambda))^{-1}\Vert)}{n}2(\Vert a_0\Vert +
\Vert \lambda \Vert + \sum_{i=1}^r \Vert a_i \Vert^2 \Vert
(Im(\lambda))^{-1}\Vert)< \frac{1}{2\Vert
(Im(\lambda))^{-1}\Vert}\right\}.$$

Since $t \mapsto \frac{P_6( t^{-1})}{n}2(\Vert a_0\Vert + t +
 \sum_{i=1}^r \Vert a_i \Vert^2 t^{-1})t^{-1}$ is a continuous
strictly decreasing function from $]0, +\infty[$ onto $]0,
+\infty[$, one can prove that ${\cal O}^{'}_n$ is an open
connected subset of $M_m(\C)$ by following the proof of (a)
Proposition 5.6 in \cite{HT}. Note that, using the inequality
$$\frac{1}{\Vert
(Im(\lambda))^{-1}\Vert} \leq \Vert \lambda \Vert,$$ one
immediately gets that for any $\lambda$ in ${\cal O}^{'}_n$,
$$\Vert (Im(\lambda))^{-1}\Vert (\Vert a_0\Vert + \Vert \lambda \Vert
+ \sum_{i=1}^r \Vert a_i \Vert^2 \Vert (Im(\lambda))^{-1}\Vert)
\geq 1$$ \noindent and thus that
$$\frac{P_6(\Vert (Im(\lambda))^{-1}\Vert)}{n}\leq \frac{1}{4}.$$
\noindent Consequently, for any $\lambda$ in ${\cal O}^{'}_n$,
$G_n(\lambda)$ is  invertible and (\ref{MastIn3}) holds. Defining
for any $\lambda$ in ${\cal O}^{'}_n$,
$$\Lambda_n(\lambda)= a_0 + \sum_{i=1}^r a_i G_n( \lambda)a_i
+ G_n (\lambda)^{-1} $$ \noindent and sticking to the proof of
\cite{HT} described in section II. 2.1, we get that, for any
$\lambda$ in ${\cal O}$,
\begin{eqnarray}
\Vert G_n(\lambda)-G(\lambda)\Vert& \leq& 4 \frac{P_6(\Vert
(Im(\lambda))^{-1}\Vert)}{n}2(\Vert a_0\Vert + \Vert \lambda \Vert
+ \sum_{i=1}^r \Vert  a_i \Vert^2 \Vert
(Im(\lambda))^{-1}\Vert)\Vert (Im(\lambda))^{-1}\Vert^2 \nonumber \\
& \leq & (\Vert \lambda \Vert+K)\frac{P_9(\Vert
(Im(\lambda))^{-1}\Vert)}{n}\label{estimdif}
\end{eqnarray}
Note that, in the following we will use (\ref{MastIn3}) in the
simplest form: $$\forall \lambda \in {\cal O}^{'}_n, ~~ \Vert
\sum_{i=1}^r a_i G_n( \lambda)a_i  \ + \ (a_0 - \lambda)+ G_n
(\lambda)^{-1}    \Vert \leq (\Vert \lambda
\Vert+K)\frac{P_7(\Vert
(Im(\lambda))^{-1}\Vert)}{n}\label{MastIn4}.$$
 \subsection{Convergence of $R_n(\lambda)$}

 \vspace{.3cm}
 Let $x_i$, $i = 1, \ldots, r$ be  self-adjoint operators  in a $C^*$ probability space $({\cal B}, \tau)$.
 We assume that the $x_i$ are free and identically semi-circular distributed
 with mean 0 and variance 1.
% Set
% $$ s= a_0 \otimes 1_{\cal B} + \sum_{p=1}^r a_p \otimes x_p \ \in M_m(C) \otimes {\cal B}, $$
% and
% $$G(\lambda) = (id_m \otimes \tau)[(\lambda \otimes 1_{\cal B}  - s)^{-1}].$$
 Then, $G$ satisfies  (\ref{eqS}).
%\begin{equation} \label{Stsc}
%a_0 + \sum_{p=1}^r a_p G(\lambda) a_p  + G(\lambda)^{-1} = \lambda
%\end{equation}
 \begin{proposition}\label{convD}
 Let $a$ be a matrix in $M_m(\C)$. Then,
 \begin{equation}
 \mathbb{E}[ \frac{1}{n} \sum_{k=1}^n (\lambda \otimes 1_n - S_n)^{-1}_{kk}\
  a \ (\lambda \otimes 1_n - S_n)^{-1}_{kk}]= GaG + O\left(\frac{P_{10}(\Vert (Im(\lambda))^{-1}\Vert)}{n}(\Vert \lambda \Vert + K)\right)
 \end{equation}
 \end{proposition}
 {\bf Proof:}
 We start from the resolvent identity:
 $$ \lambda (\lambda \otimes 1_n - S_n)^{-1}_{kk} = 1_m + \sum_{l=1}^n (S_n)_{kl} (\lambda \otimes 1_n - S_n)^{-1}_{lk}.$$
 We write $G^{(n)} (\lambda) = (\lambda \otimes 1_n - S_n)^{-1}$ and $D_a^{(n)}(\lambda) =
 \frac{1}{n} \sum_{k=1}^n G^{(n)} (\lambda) _{kk} \ a \ G^{(n)} (\lambda) _{kk}$.

 \noindent
 From the above identity,
\begin{eqnarray*}
\lefteqn{
 \lambda \ \frac{1}{n} \sum_{k=1}^n G^{(n)} (\lambda)_{kk} \ a \ G^{(n)} (\lambda)_{kk}  } \\
 && = a \ \frac{1}{n} \sum_{k=1}^n G^{(n)} (\lambda)_{kk} +
   \frac{1}{n} \sum_{k,l =1}^n (S_n)_{kl} G^{(n)} (\lambda)_{lk} \ a \ G^{(n)} (\lambda)_{kk} \\
  && =  a H_n (\lambda) + a_0 D_a^{(n)} (\lambda) + \frac{1}{n} \sum_{p=1}^r a_p \sum_{k,l =1}^n (X_n^{(p)})_{kl} G^{(n)} (\lambda)_{lk}\ a \ G^{(n)} (\lambda)_{kk}
 \end{eqnarray*}
 We take the expectation and we use the integration by part formula \eqref{IPP} for the last term:
\begin{eqnarray*}
\mathbb{E}((X_n^{(p)})_{kl} \Phi(X_n^{(1)}, \ldots, X_n^{(r)})] &= &\mathbb{E}(\Tr_n(X_n^{(p)}E_{lk}) \Phi(X_n^{(1)}, \ldots, X_n^{(r)})]\\
&=&  \frac{1}{n} \mathbb{E}[\Phi'_p(X_n^{(1)}, \ldots,
X_n^{(r)}).E_{lk}] + O\left(\frac{P_{4}(\Vert
(Im(\lambda))^{-1}\Vert)}{n^{2}}\right)
 \end{eqnarray*}
 with $\Phi(X_n^{(1)}, \ldots, X_n^{(r)}) = G^{(n)} (\lambda)_{lk}\ a \ G^{(n)} (\lambda)_{kk}$.
 Then,
 \begin{eqnarray*}
 \mathbb{E}[ (X_n^{(p)})_{kl} G^{(n)} (\lambda)_{lk}a G^{(n)} (\lambda)_{kk}] &=
 & \frac{1}{n}\mathbb{E}[G^{(n)} (\lambda)_{ll}a_p G^{(n)} (\lambda)_{kk} a G^{(n)} (\lambda)_{kk}] \\&&+
  \frac{1}{n} \mathbb{E}[G^{(n)} (\lambda)_{lk}a G^{(n)} (\lambda)_{kl} a_p G^{(n)} (\lambda)_{lk}] \\&&+ O\left(\frac{P_{4}(\Vert
(Im(\lambda))^{-1}\Vert)}{n^{2}}\right)
 \end{eqnarray*}
 Thus, we obtain from the resolvent identity,
 \begin{eqnarray*}
 \lefteqn{ (\lambda -a_0) \mathbb{E}(D_a^{(n)}(\lambda))  = a \mathbb{E}(H_n(\lambda)) + } \\
 && \sum_{p=1}^r a_p \sum_{k,l =1}^n  \frac{1}{n^2} \mathbb{E}[G^{(n)} (\lambda)_{ll}a_p G^{(n)} (\lambda)_{kk} a G^{(n)} (\lambda)_{kk} + \\
 &&  G^{(n)} (\lambda)_{lk}a G^{(n)} (\lambda)_{kl} a_p G^{(n)} (\lambda)_{lk}] + O\left(\frac{P_{4}(\Vert
(Im(\lambda))^{-1}\Vert)}{n}\right).
 \end{eqnarray*}
>From Lemma \ref{lem2}, $$\frac{1}{n^2} \sum_{k,l =1}^n \mathbb{E}
[G^{(n)} (\lambda)_{lk}a G^{(n)} (\lambda)_{kl} a_p G^{(n)}
(\lambda)_{lk}]= O\left(\frac{P_{3}(\Vert
(Im(\lambda))^{-1}\Vert)}{n}\right),$$ \noindent  thus,
$$(\lambda -a_0) \mathbb{E}(D_a^{(n)}(\lambda))  = a \mathbb{E}(H_n(\lambda)) +
  \sum_{p=1}^r a_p \mathbb{E}[H_n(\lambda) a_p D^{(n)}(\lambda)] + O\left(\frac{P_{4}(\Vert
(Im(\lambda))^{-1}\Vert)}{n}\right)$$
  From the estimate of the variance of $H_n$, we have:
  %\marginpar{ \em on utilise l'estimation de la variance}
  $$ \mathbb{E}[H_n(\lambda) a_p D_a^{(n)}(\lambda)] = \mathbb{E}[H_n(\lambda)] a_p \mathbb{E}[D_a^{(n)}(\lambda)] + O\left(\frac{P_{4}(\Vert
(Im(\lambda))^{-1}\Vert)}{n}\right).$$ Then, using also the
estimation of $\Vert G_n(\lambda)-G(\lambda)\Vert$ we get
$$(\lambda -a_0 - \sum_{p=1}^r a_p G(\lambda)
a_p)\mathbb{E}(D_a^{(n)}(\lambda))=aG +O\left(\frac{P_{9}(\Vert
(Im(\lambda))^{-1}\Vert)}{n}(\Vert \lambda \Vert +K
 )\right).$$ Using (\ref{eqS}) and (\ref{norme}) we finally get
$$\mathbb{E}(D_a^{(n)}(\lambda))=GaG + O\left(\frac{P_{10}(\Vert
(Im(\lambda))^{-1}\Vert)}{n}(\Vert \lambda \Vert +K
 )\right).$$
   $\Box$

>From the above proposition, we obtain:
\begin{proposition}
$R_n(\lambda)$ defined in Theorem \ref{ME} converges as $n$ tends to infinity to
$$R(\lambda) =\frac{\kappa_4}{2} \sum_{p=1}^r a_p G(\lambda) a_p G(\lambda) a_p G(\lambda) a_p G(\lambda).$$
More precisely, \begin{equation}\label{precision} \Vert
R_n(\lambda)- R(\lambda) \Vert \leq (\Vert \lambda \Vert +K)^2
\frac{P_{20}(\Vert (Im(\lambda))^{-1}\Vert)}{n}.\end{equation}
\end{proposition}
{\bf Proof:} It's enough to prove the convergence of each
coordinate of the $m\times m$ matrix $R_n(\lambda)$. This will
actually follow from  the convergence of terms of the form:
\begin{equation}\label{neg} \mathbb{E}\left[ _{\alpha, \beta}\!\left(n^{-1} \sum_{k=1}^n
G^{(n)}_{kk} a G^{(n)}_{kk} \right) \ _{\gamma,
\delta}\!\left(n^{-1} \sum_{k=1}^n G^{(n)}_{kk} b G^{(n)}_{kk}
\right)\right]\end{equation} for $a$ et $b$ elements of the
canonical basis in ${\cal M}_m(\C)$. Since, applying Lemma
\ref{variance}, we have
$$\mathbb{E}(\Vert D^{(n)}_a(\lambda)-\mathbb{E}(D^{(n)}_a(\lambda))\Vert^2) \leq
\frac{P_{6}(\Vert (Im(\lambda))^{-1}\Vert)}{n},$$ \noindent the
above quantity (\ref{neg}) is of the same order as:
$$ \mathbb{E}\left[ _{\alpha, \beta}\! \left(n^{-1} \sum_{k=1}^n G^{(n)}_{kk} a G^{(n)}_{kk} \right)\right]
\mathbb{E} \left[ \ _{\gamma, \delta}\!\left(n^{-1} \sum_{k=1}^n
G^{(n)}_{kk} b G^{(n)}_{kk} \right)\right] $$ According to
Proposition  \ref{convD}, this last quantity converges towards
$_{\alpha, \beta} (GaG) \ _{\gamma, \delta} (GbG)$. Thus, the
convergence of $R_n$ to $R$ follows with the estimation
(\ref{precision}). $\qquad \Box$\\

\noindent We define
$$ L(\lambda)= (id_m \otimes \tau)[(\lambda \otimes 1_{\cal B}  -
s)^{-1}(R(\lambda)G(\lambda)^{-1}\otimes 1_{\cal B})(\lambda
\otimes 1_{\cal B}  - s)^{-1}].$$

\subsection{Estimation of $ \Vert G(\lambda)-G_n(\lambda) +\frac{1}{n} L(\lambda)
\Vert$} Following 4.24 in \cite{S}, one gets for any $\lambda$ in
${\cal O}^{'}_n$,\\

$ \Vert G(\lambda)-G_n(\lambda) +\frac{1}{n} L(\lambda) \Vert$
\begin{eqnarray*}~~& \leq &\Vert (\lambda \otimes 1_{{\cal B}} - s)^{-1}\Vert
\Vert (\Lambda_n(\lambda) \otimes 1_{{\cal B}} - s)^{-1}\Vert
\Vert \Lambda_n(\lambda)-\lambda
+\frac{1}{n}R(\lambda)G(\lambda)^{-1}\Vert \\&&+ \frac{1}{n}\Vert
(\lambda \otimes 1_{{\cal B}} - s)^{-1}\Vert \Vert
R(\lambda)G(\lambda)^{-1}\Vert \Vert (\lambda \otimes 1_{{\cal B}}
- s)^{-1}- (\Lambda_n(\lambda) \otimes 1_{{\cal B}} -
s)^{-1}\Vert\\
& \leq &\Vert (Im(\lambda))^{-1}\Vert \Vert
(Im(\Lambda_n(\lambda)))^{-1}\Vert \Vert
\Lambda_n(\lambda)-\lambda
+\frac{1}{n}R(\lambda)G(\lambda)^{-1}\Vert \\&&+ \frac{C}{n}\Vert
(Im(\lambda))^{-1}\Vert^7 \Vert (Im(\Lambda_n(\lambda)))^{-1}\Vert
\Vert \Lambda_n(\lambda) - \lambda\Vert(\Vert  \lambda \Vert +K)^2
\end{eqnarray*}
where we made use of the estimates (\ref{normeG}), (\ref{Gmoins}),
(\ref{difStieljes}) and the upper bound
$$\Vert R(\lambda) \Vert \leq C\Vert
(Im(\lambda))^{-1}\Vert^4.$$ Now, for any $\lambda $ in ${\cal
O}^{'}_n$,
$$\Vert (Im(\Lambda_n(\lambda)))^{-1}\Vert \leq 2\Vert
(Im(\lambda))^{-1}\Vert$$ and
$$\Vert \Lambda_n(\lambda) - \lambda\Vert \leq \frac{P_6(\Vert
(Im(\lambda))^{-1}\Vert)}{n}(\Vert  \lambda \Vert +K).$$ Thus,
\begin{eqnarray*}
\Vert G(\lambda)-G_n(\lambda) +\frac{1}{n} L(\lambda) \Vert&\leq &
2\Vert (Im(\lambda))^{-1}\Vert^2 \Vert \Lambda_n(\lambda)-\lambda
+\frac{1}{n}R(\lambda)G(\lambda)^{-1}\Vert\\& & +
\frac{P_{14}(\Vert (Im(\lambda))^{-1}\Vert)}{n^2}(\Vert  \lambda
\Vert +K)^3.
\end{eqnarray*}
Now, for any $\lambda $ in ${\cal O}^{'}_n$,
\begin{eqnarray*}
\Vert \Lambda_n(\lambda)-\lambda
+\frac{1}{n}R(\lambda)G(\lambda)^{-1}\Vert& \leq & \Vert
\Lambda_n(\lambda)-\lambda
+\frac{1}{n}R_n(\lambda)G_n(\lambda)^{-1}\Vert\\
&&+\frac{1}{n}\Vert R_n(\lambda) (G_n(\lambda)^{-1}
-G(\lambda)^{-1})\Vert\\&& +\frac{1}{n}\Vert (R_n(\lambda)
-R(\lambda)) G(\lambda)^{-1})\Vert\\& \leq & \frac{P_7(\Vert
(Im(\lambda))^{-1}\Vert)}{n^2}(\Vert  \lambda \Vert +K)
\\
&&+\frac{P_4(\Vert (Im(\lambda))^{-1}\Vert)}{n}\Vert
G_n(\lambda)^{-1} -G(\lambda)^{-1}\Vert
\\&& +\frac{1}{n}\Vert R_n(\lambda)
-R(\lambda)\Vert (\Vert  \lambda \Vert +K)^2 \Vert
(Im(\lambda))^{-1}\Vert
\end{eqnarray*}
where we used (\ref{MastIn}), (\ref{majGinvmat}), (\ref{majR}) and
(\ref{Gmoins}). Moreover, one easily gets
\begin{eqnarray*}
\Vert G_n(\lambda)^{-1} -G(\lambda)^{-1}\Vert &=& \Vert
G_n(\lambda)^{-1}( G(\lambda) -G_n(\lambda))
G(\lambda)^{-1}\Vert\\
& \leq &\Vert G_n(\lambda)^{-1}\Vert \Vert G(\lambda)
-G_n(\lambda)) \Vert \Vert G(\lambda)^{-1}\Vert.
\end{eqnarray*}
Consequently, using the estimate (\ref{estimdif}) of $\Vert
G_n(\lambda)-G(\lambda)\Vert$ together with (\ref{majGinvmat}) and
(\ref{Gmoins}), we get
$$\Vert G_n(\lambda)^{-1} -G(\lambda)^{-1}\Vert \leq (\Vert  \lambda \Vert +K)^4 \frac{P_{11}(\Vert
(Im(\lambda))^{-1}\Vert)}{n}.$$ We conclude that
\begin{eqnarray*} \Vert
G(\lambda)-G_n(\lambda) +\frac{1}{n} L(\lambda) \Vert & \leq &
(\Vert  \lambda \Vert +K)^4 \frac{P_{17}(\Vert
(Im(\lambda))^{-1}\Vert)}{n^2}\\&&+ \frac{2}{n}\Vert R_n(\lambda)
-R(\lambda)\Vert (\Vert  \lambda \Vert +K)^2 \Vert
(Im(\lambda))^{-1}\Vert^3 .\end{eqnarray*}
%Therefore, we need an
%estimate of order $\frac{1}{n}$ of $\Vert R_n(\lambda)
%-R(\lambda)\Vert$. Actually, it is sufficient to get such an
%estimate for $\Vert E(D^{(n)}_a(\lambda) c D^{(n)}_b(\lambda)
%-GaGcGbG)\Vert$, where $a,b,c$ are any matrices in $M_m(\C)$ and
%$$D^{(n)}_a(\lambda) =
 %\frac{1}{n} \sum_{k=1}^n (\lambda \otimes 1_n - S_n)^{-1}_{kk} \ a \ (\lambda \otimes 1_n -
% S_n)^{-1}_{kk}.$$
%We have shown that
%$$E(D^{(n)}_a(\lambda))=GaG + O\left((\Vert  \lambda \Vert +K) \frac{P_{10}(\Vert
%(Im(\lambda))^{-1}\Vert)}{n}\right).$$ Moreover, applying Lemma
%\ref{variance}, we get
%$$E(\Vert D^{(n)}_a(\lambda)-E(D^{(n)}_a(\lambda))\Vert^2) \leq
%\frac{P_{6}(\Vert (Im(\lambda))^{-1}\Vert)}{n}.$$ It easily
%follows that
%$$\Vert E(D^{(n)}_a(\lambda) c D^{(n)}_b(\lambda)
%-GaGcGbG)\Vert \leq (\Vert  \lambda \Vert +K)^2 \frac{P_{20}(\Vert
%(Im(\lambda))^{-1}\Vert)}{n}.$$

Using (\ref{precision}), we can conclude that,  for any $\lambda$
in ${\cal O}^{'}_n$,
$$ \Vert G(\lambda)-G_n(\lambda) +\frac{1}{n} L(\lambda) \Vert \leq (\Vert  \lambda \Vert +K)^4 \frac{P_{23}(\Vert
(Im(\lambda))^{-1}\Vert)}{n^2}.$$ Now, for $\lambda$ in ${\cal
O}\setminus {\cal O}^{'}_n$,
\begin{eqnarray*}
1 &\leq & 4\frac{P_6(\Vert (Im(\lambda))^{-1}\Vert)}{n}(\Vert
a_0\Vert + \Vert \lambda \Vert + \Vert \sum_{i=1}^r a_i^2 \Vert
\Vert (Im(\lambda))^{-1}\Vert)\Vert
(Im(\lambda))^{-1}\Vert\}\\&\leq& (\Vert \lambda \Vert +K)
\frac{P_{8}(\Vert (Im(\lambda))^{-1}\Vert)}{n}. \end{eqnarray*} We
get
\begin{eqnarray*} \Vert
G(\lambda)-G_n(\lambda) +\frac{1}{n} L(\lambda) \Vert & \leq &
\Vert G(\lambda)-G_n(\lambda)\Vert  + \frac{1}{n}\Vert L(\lambda)
\Vert\\&\leq & (\Vert \lambda \Vert +K) \frac{P_{8}(\Vert
(Im(\lambda))^{-1}\Vert)}{n}\\ && \times \left[ (\Vert \lambda
\Vert +K) \frac{P_{9}(\Vert (Im(\lambda))^{-1}\Vert)}{n} +
\frac{1}{n}\Vert
(Im(\lambda))^{-1}\Vert^7 (\Vert \lambda \Vert +K)^2 \right]\\
&\leq & (\Vert \lambda \Vert +K)^3 \frac{P_{17}(\Vert
(Im(\lambda))^{-1}\Vert)}{n^2}.
\end{eqnarray*}
Thus, one can easily see that one can choose $K$ and $P_{23}$ such
that for any $\lambda$ in ${\cal O}$,
\begin{equation}\label{estimdifL} \Vert G(\lambda)-G_n(\lambda)
+\frac{1}{n} L(\lambda) \Vert \leq (\Vert \lambda \Vert +K)^4
\frac{P_{23}(\Vert (Im(\lambda))^{-1}\Vert)}{n^2}.\end{equation}
Note that, since under our hypothesises, $S_n$ and $-S_n$ are
identically distributed, the arguments of \cite{S} to prove her
theorem 4.5 still hold. Thus, (\ref{estimdifL}) is also valid for
any $\lambda$ such that $Im \lambda$ is negative definite.

%%%%%%%%%%%%%%%%%%%%%%

\subsection{Spectrum of $S_n$}
$\bullet$ {\bf From step 2 to step 3}\\ Sticking to the proof of
Lemma 5.5 of \cite{S}, we get that,
$$l(\lambda):= tr_m(L(\lambda1_m), \lambda \in \C \setminus \R,$$
is the Stieljes transform of a distribution $\Lambda$ with compact
support in $sp(s)$. Hence, the proof described in section 3.2
still holds (with $\alpha =4 $ and $k=23$); thus we can state that
for any smooth function $\varphi$ with compact support
\begin{equation}
\mathbb{E}[(\tr_m \otimes \tr_n ) (\varphi(S_n))] = (\tr_m \otimes
\tau) (\varphi(s))+ \frac{1}{n}\Lambda(\varphi) +
O(\frac{1}{n^2}).
\end{equation}
Moreover, following the proof of Lemma 5.6 in \cite{S}, one can
show that $\Lambda(1)=0$ and deduce that, for $\varphi$ smooth,
constant outside a compact set and such that $supp(\varphi) \cap
sp(s) = \emptyset$,
$$\mathbb{E}[(\tr_m \otimes \tr_n ) (\varphi(S_n))] =O(\frac{1}{n^2}).$$
$\bullet$ {\bf Step 4}\\
The proof of step 4 is exactly the same as in \cite{HT} so that we
have proved that,
 for any $\varepsilon >0$ and
almost surely
$$Spect(S_n) \subset Spect(s)+(-\varepsilon,\varepsilon)$$
when $n$ goes to infinity. Note that this result implies that
$$sup_n \Vert X_n^{(p)}\Vert < + \infty ~~a.e.$$

%%%%%%%%%%%%%%%%%%%%%%%%%%%%%%%%
\subsection{The main theorem}
\subsubsection{First inequality}
By the same arguments of \cite{HT} in Proposition 7.3, we can
deduce the following inequality from the above inclusion of the
spectrum of $S_n$.
\begin{proposition}
Almost everywhere,
for all polynomials $p$ in $r$ non commuting variables,
$$ \limsup_{n \rightarrow + \infty} \Vert
p(X_n^{(1)},\ldots,X_n^{(r)}\Vert \leq \Vert
p(x_1,\ldots,x_r)\Vert.$$
\end{proposition}
\subsubsection{Second inequality}
\begin{proposition}
Almost everywhere,
for all polynomials $p$ in $r$ non commuting variables,
$$ \liminf_{n \rightarrow + \infty} \Vert
p(X_n^{(1)},\ldots,X_n^{(r)}\Vert \geq \Vert
p(x_1,\ldots,x_r)\Vert.$$
\end{proposition}
>From \cite{HT}, Proof of Lemma 7.2, it is clear that this
proposition follows from the almost sure asymptotic freeness of
the $X_n^{(i)}$ together with the property that
$$a.e \sup_n \Vert X_n^{(i)} \Vert < +\infty.$$
The proof of the first point follows the proof of Theorem 6.2 in
\cite{S}; nevertheless, we modify the proof of Lemma 6.5 in
\cite{S} to get the analogue in our context without needing such a
result as Lemma 6.4 in \cite{S}.
\begin{lemme}\label{poincaretrace}
Let $d$ be in $\N^{*}$, $i_1, \ldots,i_d$ be in $\{1, \ldots,r\}$ and $n$ be in $\N^{*}$.
Define $f: M_n(\C)^r \rightarrow \C$ by
$$f(v_1, \ldots, v_r) = \tr_n(v_{i_1} \ldots v_{i_d}).$$
\noindent
Then, there is a constant $C>0$ such that
$$\mathbf{V}{f(X_n^{(1)}, \ldots,
X_n^{(r)})} \leq \frac{C}{n^2}.$$
\end{lemme}
{\bf Proof:}
Applying Poincar\'e Inequality (\ref{Poincare}), we get
$$\mathbf{V}{f(X_n^{(1)}, \ldots,
X_n^{(r)})} \leq \frac{C}{n}\mathbb{E}\{\Vert \nabla f(X_n^{(1)},
\ldots, X_n^{(r)}) \Vert_e^2 \}.$$ \noindent Now, let $v=(v_1,
\ldots, v_r)$ be in $M_n(\C)^r $ and $w=(w_1, \ldots, w_r)$ be in
$M_n(\C)^r_{sa} $ with $\Vert w \Vert_e =1$. By the Cauchy
Schwartz inequality
\begin{eqnarray*}
\vert {\frac{d~}{dt}}_{t=0} f(v+tw) \vert & =& \frac{1}{n}\vert
\Tr_n(w_{i_1}v_{i_2}v_{i_3}\ldots v_{i_d})
+\Tr_n(v_{i_1}w_{i_2}v_{i_3}\ldots v_{i_d})+ \ldots \\ & & ~~~~~~~~~~~~~~~~~~~~~~~~~~~~~~~~~~~~~~~~~~\ldots   +\Tr_n(v_{i_1}v_{i_2} \ldots v_{i_{d-1}} w_{i_d}) \vert\\
& \leq & \frac{1}{n}\left\{  \Vert w_{i_1}\Vert_{2}  \Vert
v_{i_2}v_{i_3}\ldots v_{i_d}\Vert_{2}
%+\Vert w_{i_2}\Vert_{2}    \Vert v_{i_1}v_{i_3}\ldots v_{i_d}\Vert_{2}
 + \ldots    +
\Vert w_{i_d}\Vert_{2}    \Vert v_{i_1}v_{i_2} \ldots v_{i_{d-1}}\Vert_{2}\right\}  \\
& \leq & \frac{1}{n}\left\{  \Vert  v_{i_2}v_{i_3}\ldots
v_{i_d}\Vert_{2} + \Vert v_{i_1}v_{i_3}\ldots v_{i_d}\Vert_{2}  +
\ldots    +
    \Vert v_{i_1}v_{i_2} \ldots v_{i_{d-1}}\Vert_{2}\right\}
\end{eqnarray*}
Thus,
\begin{eqnarray*}
\Vert\nabla f(X_n^{(1)}, \ldots, X_n^{(r)}) \Vert_e^2  &\leq &
\frac{1}{n^2}\left[ \left\{  \Tr_n(X_n^{(i_2)}X_n^{(i_3)}\ldots X_n^{(i_d)}X_n^{(i_d)}\ldots X_n^{(i_2)})\right\}^{\frac{1}{2}} \right.\\
& &\left. ~~~~~~+ \ldots + \left\{  \Tr_n(X_n^{(i_1)}\ldots X_n^{(i_{d-1})}X_n^{(i_{d-1})}\ldots X_n^{(i_1)})\right\}^{\frac{1}{2}}\right]^2\\
&\leq & \frac{C}{n^2} \left\{  \Tr_n(X_n^{(i_2)}X_n^{(i_3)}\ldots X_n^{(i_d)}X_n^{(i_d)}\ldots X_n^{(i_2)})\right.\\
& & \left.~~~~~~+ \ldots +   \Tr_n(X_n^{(i_1)}\ldots
X_n^{(i_{d-1})}X_n^{(i_{d-1})}\ldots X_n^{(i_1)})\right\}
\end{eqnarray*}
\noindent for some constant $C$ depending on $d$, and we get that
\begin{eqnarray*}
 \mathbb{E} \left(\Vert \nabla f(X_n^{(1)},
\ldots, X_n^{(r)}) \Vert_e^2 \right)& \leq & \frac{C}{n} \left\{
\mathbb{E} \left( \tr_n(X_n^{(i_2)}X_n^{(i_3)}\ldots X_n^{(i_d)}X_n^{(i_d)}\ldots X_n^{(i_2)})\right)\right.\\
& &\left. ~~~~~~+ \ldots +  \mathbb{E} \left(
\tr_n(X_n^{(i_1)}\ldots X_n^{(i_{d-1})}X_n^{(i_{d-1})}\ldots
X_n^{(i_1)})\right) \right\}.
\end{eqnarray*}

\noindent Each term inside the brackets of the left hand side  is
uniformly bounded in $n$ since it converges as $n$ tends to
infinity according to the result of asymptotic freeness in mean of
Dykema in \cite{D}. The result follows. $\Box$

Lemma \ref{poincaretrace} yields the almost sure asymptotic
freeness of the $X_n^{(i)}$ using the Borel Cantelli lemma.

~~

\noindent In conclusion,
\begin{theoreme}
Let $X_n^{(1)},\ldots,X_n^{(r)}$ be independent Wigner matrices
associated to a symmetric distribution $\mu$ which satisfies a
Poincar\'e inequality. Let $(x_1, \ldots, x_r)$ be a semicircular
system. Then, almost everywhere, for all polynomials $p$ in $r$ non commuting
variables
$$ \lim_{n \rightarrow + \infty} \Vert
p(X_n^{(1)},\ldots,X_n^{(r)})\Vert = \Vert
p(x_1,\ldots,x_r)\Vert.$$
\end{theoreme}

%%%%%%%%%%%%%%%%%%%%%%%%%%%%%
\section{The Wishart case}
We consider a $n \times n$ Hermitian matrix $Y$, distributed as a
Wishart matrix of parameter $p(n) \geq n$ and variance
$\frac{1}{n}$ that is with density w.r.t the Lebesgue measure $dM$
on ${\cal M}_{sa} (\C)$:
$$ C_{n,p}  1_{ (M \geq 0)} (\det(M))^{p(n) -n} \exp( - n \tr(M)).$$
We assume that $\displaystyle \frac{p(n)}{n}\vers_{n\rightarrow
\infty} \alpha$ for some $\alpha \geq 1$. More precisely,
according to Dirichlet theorem (\cite{T}, Lemme 14.1), there
exists  subsequences $p(n)$ and $q(n)$ of integers tending to
$\infty$ such that:
 $$|\frac{p(n)}{q(n)} - \alpha | \leq \frac{1}{q(n)^2}.$$
 So, we shall consider a matrix $Y$ of size $q(n)$ and parameter $p(n)$. For simplicity, we shall denote
 the subsequence $q(n)$ by $n$ and therefore, we will assume in this section that:
 \begin{equation} \label{dirichlet}
 |\frac{p(n)}{n} - \alpha | \leq \frac{1}{n^2}.
 \end{equation}
It is well know that the spectral measure of $Y$ converges to the
so called Marchenko-Pastur distribution  $\mu_\alpha$ \cite{MP}:
$$ \mu_\alpha (dx) = \frac{\sqrt{((\sqrt{\alpha}+1)^2 -x)(x -(\sqrt{\alpha}-1)^2 )}}{2\pi x} 1_{[(\sqrt{\alpha}-1)^2, (\sqrt{\alpha}+1)^2]} (x) dx.$$

\subsection{ Differentiation formula for the Wishart ensemble}
\begin{lemme}
Let $\Phi$ a $C^1$ function on ${\cal M}_{sa} (\C)$ with $\Phi(0)
= 0$, then:
\begin{equation} \label{df}
\mathbb{E}\left[ \Phi'(Y).H \right] - n \mathbb{E}\left[ \Phi(Y)
\Tr_n(H)\right] + (p(n) -n)  \mathbb{E}\left[ \Phi(Y)
\Tr_n(Y^{-1}H)\right] =  0
\end{equation}
for all hermitian matrix $H$, or by linearity for $H = E_{jk}$,
$1\leq j,k \leq n$.
\end{lemme}
{\bf Proof:} Since the Lebesgue measure is invariant by
translation,
$$\mathbb{E} [\Phi(Y)] = \int \Phi(M+\epsilon X) \exp(- n \Tr_n(M+\epsilon X)) (\det(M+\epsilon X))^{p(n) - n}
1_{(M+\epsilon X \geq 0)} dM.$$ Now, by differentiation with
respect to $\epsilon$ and taking $\epsilon =0$, we obtain
\eqref{df} using $\partial(\det M)^s = s(\det M)^s M^{-1}$.

\subsection{The master equation}
Let $(X^{(1)}_n, \ldots, X^{(r)}_n)_{i = 1, \ldots r}$ be $r$
independent copies of the random matrix $Y$. We shall apply
\eqref{df} with
$$\Phi(X_n^{(l)}) = \left[ (1_m \otimes X_n^{(l)}) (\lambda \otimes 1_n - S_n)^{-1}\right]_{jk} \in {\cal M}_{m} (\C)$$
 and $ H = E_{jk}$. Then,
$$ \Phi'(X_n^{(l)}).E_{jk}  = (\lambda \otimes 1_n - S_n)^{-1}_{kk} +[ (1_m \otimes X_n^{(l)}) (\lambda \otimes 1_n - S_n)^{-1}]_{jj} a_l  [(\lambda \otimes 1_n - S_n)^{-1}]_{kk} $$
and
$$\Phi(X_n^{(l)}) \Tr_n((X_n^{(l)})^{-1}E_{jk}) = (X_n^{(l)})^{-1}_{kj}  \left[ (1_m \otimes X_n^{(l)}) (\lambda \otimes 1_n - S_n)^{-1}\right]_{jk}. $$
The sum over $j$ of the terms in the above equation gives:
$$(\lambda \otimes 1_n - S_n)^{-1}_{kk}.$$
Now, if we sum the identities obtained by \eqref{df} over $j,k$,
and dividing by $n^2$, we obtain:
\begin{eqnarray}
\lefteqn{
\mathbb{E}[ id \otimes \tr_n ((\lambda \otimes 1_n - S_n)^{-1})]  } \\
&&
 + \mathbb{E}[ id \otimes \tr_n ((1_m \otimes X_n^{(l)}) (\lambda \otimes 1_n - S_n)^{-1})\  a_l \  id \otimes \tr_n ((\lambda \otimes 1_n - S_n)^{-1})]  \nonumber\\
&& - \mathbb{E}[id \otimes \tr_n ((1_m \otimes X_n^{(l)}) (\lambda \otimes 1_n - S_n)^{-1})]  \nonumber \\
&& + (\frac{p(n)}{n} - 1) \mathbb{E}[id \otimes \tr_n ( (\lambda
\otimes 1_n - S_n)^{-1})]  =0 \nonumber
\end{eqnarray}
which can be written as:
$$ \mathbb{E} \left[ id \otimes \tr_n ((1_m \otimes X_n^{(l)}) (\lambda \otimes 1_n - S_n)^{-1}) (1_m - a_l \
id \otimes \tr_n ( (\lambda \otimes 1_n - S_n)^{-1})\right] $$
\begin{equation} \label{ippW}
= \frac{p(n)}{n} \mathbb{E}[id \otimes \tr_n( (\lambda \otimes 1_n
- S_n)^{-1})].
\end{equation}

\begin{proposition} \label{lemW}
\begin{enumerate}
\item For $\lambda \in {\cal O}$,
$$| \mathbb{E} \left[ id \otimes \tr_n ((1_m \otimes X_n^{(l)}) (\lambda \otimes 1_n - S_n)^{-1}) (1_m - a_l \
id \otimes \tr_n ( (\lambda \otimes 1_n - S_n)^{-1})\right] -  $$
$$\mathbb{E} \left[ id \otimes \tr_n ((1_m \otimes X_n^{(l)}) (\lambda \otimes 1_n - S_n)^{-1})\right] \mathbb{E}\left[ (1_m - a_l \
id \otimes \tr_n ( (\lambda \otimes 1_n - S_n)^{-1})\right] | $$ $
\leq \frac{P_4(||(Im(\lambda))^{-1}||)}{n^2}$. \item For $\lambda
\in {\cal O}$,
%$-Im(G_n(\lambda)) \geq \frac{1}{c(\lambda)} 1_m$, thus
$$ || G_n(\lambda)^{-1} || \leq (|| \lambda|| +K)^2||(Im(\lambda))^{-1}||$$
\item If $a_l$ is invertible and $\lambda \in {\cal O}$, then
$(1_m - a_l G_n(\lambda))$ is invertible and
$$ || (1_m - a_l G_n(\lambda))^{-1} || \leq ||a_l^{-1}|| (|| \lambda|| +K)^2||(Im(\lambda))^{-1}||.$$
If $||(Im(\lambda))^{-1}|| < \frac{1}{2||a_l||}$, then $(1_m - a_l
G_n(\lambda))$ is invertible and
\begin{equation} \label{majcste}
 || (1_m - a_l G_n(\lambda))^{-1} || \leq 2.
 \end{equation}
\end{enumerate}
\end{proposition}
{\bf Sketch of Proof:} \\
1. The variance estimate follows from the Gaussian Poincar\'e inequality since we can write $Y = \frac{1}{n} X^*X$ for a rectangular Gaussian matrix $X$.
 We proceed as in \cite[Section 4]{HT}. We need some estimate on the maximal eigenvalue of $Y$, i.e. $\mathbb{E}[\lambda_{max}] $
and  $\mathbb{E}[\lambda^3_{max}] $ are bounded, independently of $n$. This can be proved, as in Lemma 5.1 of \cite{HT}, using previous results in \cite{HT0}. \\
2. The proof  is the same as  Proposition 5.2 in \cite{HT}. \\
3.  If $a_l$ is invertible,
$$(1-a_l G_n(\lambda)) = a_l (a_l^{-1} - G_n(\lambda)) := a_l T.$$
Now, the matrix $T$ satisfies, $Im(T) = - Im(G_n(\lambda))$ and
thus is positive definite (see the proof of Proposition 5.2 in
\cite{HT}). Its inverse $T^{-1}$ satisfies:
$$ ||T^{-1}|| \leq ||Im(T)^{-1} || =  ||Im(G_n(\lambda))^{-1} || \leq (|| \lambda|| +K)^2||(Im(\lambda))^{-1}||.$$
\eqref{majcste} follows from the majoration $||G_n(\lambda)|| \leq
||(Im(\lambda))^{-1}||$. $\Box$

\vspace{.3cm} \noindent From \eqref{ippW} and Proposition
\ref{lemW}, we obtain:
\begin{eqnarray*}
\lefteqn{
\left \Vert \sum_{l=1}^r a_l \mathbb{E}[id \otimes \tr_n ((1_m \otimes Y^{(l)}) (\lambda \otimes 1_n - S_n)^{-1}) ] -\right.} \\
&& \left.  \frac{p(n)}{n} \sum_{l=1}^r a_l
\mathbb{E}[H_n(\lambda)] (1_m - a_l
\mathbb{E}[H_n(\lambda)])^{-1}\right\Vert \leq   (|| \lambda||
+K)^2 \frac{P_5(||(Im(\lambda))^{-1}||)}{n^2}
\end{eqnarray*}
The first line  of the above equation equals:
$$ \mathbb{E}[id \otimes \tr_n (S_n - a_0\otimes 1_n) (\lambda \otimes 1_n - S_n)^{-1}) ] =
- 1_m + (\lambda -a_0) \mathbb{E}[H_n(\lambda)].$$ We have thus
obtain the master inequality:
\begin{equation}
|| -a_0 G_n(\lambda) + \lambda G_n(\lambda)  - \frac{p(n)}{n}
\sum_{l=1}^r a_l G_n(\lambda) (1_m - a_l G_n(\lambda))^{-1}  -1_m
|| \leq (|| \lambda|| +K)^2
\frac{P_5(||(Im(\lambda))^{-1}||)}{n^2}
\end{equation}
or since the matrices $(1_m - a_l G_n(\lambda))^{-1}$ and $ a_l
G_n(\lambda)$ commute,
\begin{equation}\label{mastertineqW}
|| -a_0 G_n(\lambda) + \lambda G_n(\lambda)  - \frac{p(n)}{n}
\sum_{l=1}^r (1_m - a_l G_n(\lambda))^{-1}  a_l G_n(\lambda) \
-1_m || \leq (|| \lambda|| +K)^2
\frac{P_5(||(Im(\lambda))^{-1}||)}{n^2}
\end{equation}
%%%%%%%%%%%%%%%%%%%%%%%%%%%%%%%%%
\subsection{Estimation of $||G_n(\lambda)- G(\lambda)||$}
Let $x_i, i\leq r$ be a free family of self adjoint variables in a
$C^*$-probability space $({\cal B}, \tau)$, with Marchenko-Pastur
distribution $\mu_\alpha$, with parameter $\alpha$.

%\noindent We define
$$s= a_0 \otimes 1_{{\cal B}} + \sum_{i=1}^r a_i \otimes x_i$$
%and
%$$ G(\lambda) = (id_m \otimes \tau)[(\lambda \otimes 1_n - S_n)^{-1}] \in M_m(\C).$$
Using the known expression of the $R$ transform of the
distribution of $x_i$ (see \cite{CC}, \cite[Example
3.3.5]{HP}\footnote{We warm the reader that the $R$ transform
defined in this book differs by a factor $z$ from the Voiculescu
$R$ transform we used here}):
$$ R_x(z) = \alpha (1-z)^{-1}, z \in \C \backslash \R,$$
we can show the following
\begin{lemme}
$G$ satisfies the following equation: for $\lambda \in {\cal O}$,
\begin{equation}\label{GMP}
a_0 + \alpha \sum_{i=1}^r (1_m - a_i G(\lambda))^{-1} a_i +
G(\lambda)^{-1} = \lambda.
\end{equation}
\end{lemme}
{\bf Sketch of Proof:} From the definition of the $R$
transformation with amalgation over $M_m(\C)$, we can show that:
$$ R_{a \otimes x} (\lambda) = R_x (a \lambda) a, \; a \in M_m(\C)_{sa}, \lambda \in {\cal O}$$
and then, by freeness asumption
$$R_{ \sum_{i=1}^r a_i \otimes x_i} (\lambda) =  \sum_{i=1}^r R_{x_i}(a_i \lambda) a_i.$$
\eqref{GMP} follows, using the relation between $R$ and $G$.
$\Box$.
\begin{theoreme} \label{theoGWis}
For  any $\lambda \in {\cal O}$,
\begin{equation}
||G(\lambda) - G_n(\lambda)|| \leq (|| \lambda|| +K)^4
\frac{P_8(||(Im(\lambda))^{-1}||)}{n^2}.
\end{equation}
\end{theoreme}
{\bf Proof:} We can proceed as in the proof of Theorem 5.7 in
\cite{HT}. We  just mention the different steps:

\vspace{.3cm} \noindent {\bf Step 1}:
 Define $\Lambda_n(\lambda)  = a_0 + G_n(\lambda)^{-1} + \alpha \sum_l (1_m - a_l G_n(\lambda))^{-1} a_l$. From the master  inequality \eqref{mastertineqW}, Proposition  \ref{lemW} and \eqref{dirichlet}, we can show that:
 $$|| \lambda - \Lambda_n(\lambda) || \leq (|| \lambda|| +K)^4 \frac{P_6(||(Im(\lambda))^{-1}||)}{n^2}.$$
Then, for $\lambda \in {\cal O}'_n$ of the form
$${\cal O}'_n =\left\{ \lambda \in {\cal O}, (|| \lambda|| +K)^4 \frac{P_6(||(Im(\lambda))^{-1}||)}{n^2} < \frac{1}{2}
\frac{1}{||(Im(\lambda))^{-1}||}\right\},$$ we have
$Im(\Lambda_n(\lambda)) \geq \frac{1}{2||(Im(\lambda))^{-1}||}$
and in particular $\Lambda_n(\lambda)  \in {\cal O}$.

\vspace{.3cm} \noindent {\bf Step 2}: For $\lambda \in {\cal
O}'_n$, we can consider $G(\Lambda_n(\lambda))$ and we have, from
the identity \eqref{GMP}:
\begin{equation} \label{1.10}
 a_0 + G_n(\lambda)^{-1} + \alpha \sum_l (1_m - a_l G_n(\lambda))^{-1} a_l =
 a_0 + G(\Lambda_n(\lambda))^{-1} + \alpha \sum_l (1_m - a_l G(\Lambda_n(\lambda)))^{-1} a_l.
 \end{equation}
 \begin{lemme} (see \cite{HT}, Propostion 5.6) For $\lambda \in  {\cal O}'_n$,
 \begin{equation} \label{egW}
 G(\Lambda_n(\lambda)) = G_n (\lambda).
 \end{equation}
  \end{lemme}
 {\bf Proof:}  As in \cite{HT}, it's enough to prove \eqref{egW} for $\lambda \in  {\cal O}''_n$, a non empty subset of the connected subset $ {\cal O}'_n$.
 Put $x = G_n(\lambda)$ and $y= G(\Lambda_n(\lambda))$, then, from \eqref{1.10},
 $$a_0 + x^{-1} + \alpha \sum_l (1_m - a_l x)^{-1} a_l =
 a_0 + y^{-1} + \alpha \sum_l (1_m - a_l y)^{-1} a_l $$
 so that
 $$y + \alpha \sum_l x (1_m - a_l x)^{-1} a_l\ y =
 x + \alpha \sum_l x (1_m - a_l y)^{-1} a_l\ y,$$
Thus,
\begin{eqnarray*}
y-x &= & \alpha \sum_l x [(1_m - a_l y)^{-1}- (1_m - a_l x)^{-1}] a_l\ y \\
&=&  \alpha \sum_l x (1_m - a_l y)^{-1} [(1 -a_l x) - (1-a_l y)](1_m - a_l x)^{-1} a_l \ y \\
&=&  \alpha \sum_l x (1_m - a_l y)^{-1}a_l (y-x) (1_m - a_l
x)^{-1} a_l\  y
\end{eqnarray*}
In particular, we have,
\begin{equation} \label{contraction}
 \Vert y-x \Vert \leq \left(  \alpha \Vert x \Vert \Vert y \Vert  \sum_l   \Vert(1_m - a_l y)^{-1} \Vert \Vert(1_m - a_l x)^{-1} \Vert \Vert a_l \Vert^2 \right) \Vert y-x \Vert
 \end{equation}
 Now, we have
 $$\Vert x \Vert = \Vert G_n(\lambda) \Vert \leq ||(Im(\lambda))^{-1}||$$ and
$$\Vert y \Vert = \Vert G(\Lambda_n(\lambda)) \Vert \leq (||(Im(\Lambda_n(\lambda)))^{-1}|| \leq
2 ||(Im(\lambda))^{-1}|| $$
for $\lambda \in {\cal O}'_n$ (see Step 1). \\
Moreover, from Proposition \ref{lemW}, for
$||(Im(\lambda))^{-1}||$ small enough,
$$  \Vert(1_m - a_l x)^{-1} \Vert  \leq 2; \;  \Vert(1_m - a_l y)^{-1} \Vert \leq 2.$$
Set
$$ {\cal O}''_n = \{ \lambda \in {\cal O}'_n;  \alpha \Vert G_n(\lambda)  \Vert
\Vert G(\Lambda_n(\lambda)) \Vert  \sum_l   \Vert(1_m - a_l
G(\Lambda_n(\lambda)))^{-1} \Vert \Vert(1_m - a_l
G_n(\lambda))^{-1} \Vert \Vert a_l \Vert^2 <1 \},$$ Then, from
\eqref{contraction}, for $\lambda \in {\cal O}''_n$,
$G(\Lambda_n(\lambda)) = G_n(\lambda)$. Now, it is easy to see,
from the above estimates,  that $\lambda = it1_m \in {\cal O}''_n$
for $t$ large enough, so  ${\cal O}''_n $ is a non empty set.
$\Box$

 %%%%%%%%%%%%%%%%%%%%%%%%%
%%%%%%%%%%%%%%%%%%%%%%%%%
\vspace{.3cm} \noindent {\bf Step 3}:
 The estimation of $G(\Lambda_n(\lambda)) - G(\lambda)$ is obtained as in Subsection 3.1 (considering the two cases $\lambda \in {\cal O}'_n$ and
 $\lambda \in {\cal O} \backslash {\cal O}'_n$).  $\Box$
%%%%%%%%%%%%%%%%%%%%%
\subsection{The spectrum of $S_n$}
>From Theorem \ref{theoGWis} and the proof described in Section 3.2
(see also Section 6 in \cite{HT}), we can prove that for $\varphi$
smooth, constant outside a compact set and such that
$supp(\varphi) \cap sp(s) = \emptyset$
$$\mathbb{E}[(\tr_m \otimes \tr_n ) (\varphi(S_n))] =O(\frac{1}{n^2}).$$
from which we deduce that,  for any $\varepsilon >0$ and almost
surely
$$Spect(S_n) \subset Spect(s)+(-\varepsilon,\varepsilon)$$
when $n$ goes to infinity.
%%%%%%%%%%%%%%%%%
\subsection{The main theorem}
We can now prove:
\begin{theoreme}
There exists a set $N$ of probability 0 such that for all non
commutative polynomial $p$ in $r$ variables, and all $\omega \in
\Omega \backslash N$,
\begin{equation} \lim_{n \vers \infty} || p(X_n^{(1)}(\omega), \ldots, X_n^{(r)}(\omega)) || = || p(x_1, \ldots x_r)||.
\end{equation}
\end{theoreme}
{\bf Proof:} The inequality
\begin{equation} \label{5.14}
\limsup_{n \vers  \infty} || p(X_n^{(1)}, \ldots, X_n^{(r)}) ||
\leq  || p(x_1, \ldots x_r)|| \ a.s.
\end{equation}
follows from the above inclusion of the spectrum of $S_n$ and the
arguments developed in \cite{HT}, Section 7. The reverse
inequality
\begin{equation}
\liminf_{n \vers  \infty} || p(X_n^{(1)}, \ldots, X_n^{(r)}) ||
\geq  || p(x_1, \ldots x_r)|| \ a.s.
\end{equation}
follows, as in Lemma 7.2 in \cite{HT},  from the a.s. asymptotic
freeness of the $(X_n^{(i)})_{i= 1,\ldots,r}$ and $\sup_n \Vert
X_n^{(i)} \Vert < \infty$ a.s.. The first point was proved by Hiai
and Petz (see \cite{HP}, \cite{HP1}) and the second point follows
from \eqref{5.14}. $\Box$

\vspace{.5cm} \noindent {\bf Remark:} If we only assume the
convergence of $\frac{p(n)}{n}$ to $\alpha$ with $|\frac{p(n)}{n}
- \alpha | \leq \frac{1}{n}$, then an extra term appears in the
estimation of $G-G_n$ at order $n^{-2}$, namely:
$$ || G(\lambda) - G_n(\lambda) + (\alpha- \frac{p(n)}{n}) L(\lambda) || \leq \frac{C(\lambda)}{n^2}$$
with
$$ L(\lambda) = (id_m \otimes \tau)[(\lambda \otimes 1_{{\cal B}} -s)^{-1} (R(\lambda) \otimes 1_{{\cal B}})
(\lambda \otimes 1_{{\cal B}} -s)^{-1}]$$
and $R(\lambda) = \sum_l (1_m - a_l G(\lambda))^{-1} a_l $. \\
As in Schultz \cite{S} and in the iid case (see Section 4), this
term gives rise to a distribution with compact support in $sp(s)$
and the conclusion remains true.


\begin{thebibliography}{99}
\bibitem{Tou}
%MR1845806 (2002g:46132) 46N20 (26D15 46E99 58J60 60J10 60J60)
An\'e, C., Blach\`ere, S., Chafa\"\i, D., Foug\`eres, P., Gentil,
I., Malrieu, F., Roberto, C., Scheffer, G. {\it Sur les
in\'egalit\'es de Sobolev logarithmiques.} (French) [Logarithmic
Sobolev inequalities] Panoramas et Synth\`eses [Panoramas and
Syntheses], 10. Soci\'et\'e Math\'ematique de France, Paris, 2000.
%xvi+217 pp. ISBN: 2-85629-105-8
\bibitem{Ba} Bai, Z. D.: {\it Methodology in Spectral analysis of large
dimensional random matrices. A review.} Statistica Sinica {\bf 9}
(1999), pp 611-677.
\bibitem{BY} Bai, Z.D., Yin, Y. Q.: {\it Necessary and sufficient
conditions for almost sure convergence of the largest eigenvalue
of a Wigner matrix.} Ann. of Proba. {\bf 16} (1988), pp 1729-1741.
\bibitem{B} Bobkov, S., Gotze.: {\it Exponential integrability and transportation cost related to logarithmic Sobolev inequalities.} J. Funct. Anal.  {\bf 163} (1999), pp 1-28.
\bibitem{CC} Capitaine, M. and  Casalis, M.:{\it Asymptotic freeness by generalized moments for Gaussian and Wishart matrices.
Application to beta random matrices.}
 Indiana Univ. Math. J. {\bf  53} (2004), pp 397--431.
\bibitem{D} Dykema, K.: {\it  On Certain Free Product Factors via an Extended Matrix Model.}
 J. Funct. Anal. {\bf 112} (1993), pp 31-60.
 \bibitem{G} Geman S.: {\it A limit theorem for the norm of random
 matrices.} Annals of Probability 8 (1980) pp 252-261.
\bibitem{HT0} Haagerup, U. and Thorbj\o rnsen, S.: {\it Random matrices with complex Gaussian
entries.}
 Expo. Math. {\bf  21} (2003), pp 293--337.
\bibitem{HT} Haagerup, U. and  Thorbj\o rnsen, S.: {\it A new application of random matrices: $Ext(C^*_{red}(F_2))$ is not a group.} To appear in Ann. Math.
\bibitem{HP} Hiai, F. and Petz, D.: "The semicircle law, free random variables and entropy." Mathematical Surveys and Monographs, 77. American Mathematical Society, Providence, RI, 2000.
\bibitem{HP1} Hiai, F. and Petz, D.: {\it Asymptotic freeness almost everywhere for random matrices.} Acta Sci. Math.  {\bf  66} (2000)  pp 809--834.
\bibitem{KKP} Khorunzhy, A., Khoruzhenko, B., Pastur,L.:
{\it Asymptotic properties of large random matrices with
independent entries}, J. Math. Phys. {\bf  37} (1996),
pp 5033-5060.
\bibitem{L} Ledoux, L.: "The Concentration of Measure Phenomenon."
Mathematical Surveys and Monographs, Volume 89, A.M.S, 2001.
\bibitem{MP} Marchenko, V. and Pastur, L.: {\it The distribution of eigenvalues in a certain sets of random matrices}.  Math. Sb.  {\bf 72} (1967), pp 507-536.
\bibitem{S} Schultz, H.: {\it Non-commutative polynomials of independent Gaussian  random matrices. The real and symplectic cases.} To appear in {\it Prob. Th. Rel. Fields}.
\bibitem{T} Tenenbaum, G.: "Introduction \`a la th\'eorie analytique et probabiliste des nombres."  Institut Elie Cartan, Nancy, 1990.
\bibitem{Th} Thorbj\o rnsen, S.:{\it Mixed moments of Voiculescu's Gaussian Random matrices.} J. Funct. Anal. {\bf 176} (2000), pp 213-246.
\bibitem{V} Voiculescu, D.:  {\it Limit laws for random matrices and free products.} Invent. Math. {\bf 104} (1991), pp 201-220.
\end{thebibliography}
\end{document}